\documentclass[12pt]{amsart}
\usepackage{amssymb}
\usepackage[usenames]{color}

\newtheorem{theorem}{Theorem}

\newtheorem{proposition}[theorem]{Proposition}

\newtheorem{corollary}[theorem]{Corollary}
\newtheorem{remark}[theorem]{Remark}
\newtheorem{definition}{Definition}
\newtheorem{lemma}[theorem]{Lemma}

\newcommand{\G}{{\mathfrak G}}

\begin{document}

{

\title[Isotopy classification of degree 3 Morse polynomials in ${\mathbb R}^3$]{Isotopy classification of Morse polynomials of degree three on ${\mathbb R}^3$}

\author{V.A.~Vassiliev}
\address{Weizmann Institute of Science, Rehovot, Israel}
\thanks{This work was supported by the Absorption Center in Science of the Ministry of Immigration and Absorption of the State of Israel}
 \email{vavassiliev@gmail.com}

\begin{abstract}
We enumerate all isotopy classes of degree three Morse polynomials ${\mathbb R}^3 \to {\mathbb R}^1$ with nonsingular principal homogeneous parts, proving that there are exactly 37 of them. We also count all 2258 isotopy classes of {\em strictly} Morse polynomials ${\mathbb R}^3 \to {\mathbb R}^1$ of degree three with the maximal possible number (eight) of real critical points. 
A main tool in  this classification is a combinatorial computer program that formalizes Morse surgeries, local monodromy and Picard-Lefschetz theory. 
\end{abstract}


\maketitle

\section{Introduction}

\subsection{Main results}
Two of the main objects of bifurcation theory are {\em discriminants} and {\em caustics}, which are sets of functions with a singular zero level set and with non-Morse critical points, respectively. These objects are usually considered in parallel and in close relation to each other, see, for example, \cite{AVGZ1}, \cite{Kluwer}. Much of this study is devoted to enumerating the connected components of their complementary sets of functions without such degeneracies. This problem for discriminants, i.e., the classification of real polynomials with smooth zero-level sets, has been very popular since Hilbert's 16th problem, see e.g.  \cite{klein}, \cite{Pe}, \cite{Ro}, \cite{DK}, \cite{FK}, \cite{Kha}, \cite{schlaf}, \cite{viro1}, \cite{viro2}, \cite{zeuten}. We study the analogous problem for caustics: the enumeration of components of spaces of Morse polynomials. 

For an earlier step in this direction, see the papers \cite{ANN}, \cite{AN}, which contain some estimates of the number of possible topological types of {\em strictly} Morse polynomials ${\mathbb R}^2 \to {\mathbb R}^1$ of degree four.  In the corresponding local theory, all topological types of arbitrarily small Morse perturbations of real simple singularities were enumerated, see \cite{gor}. 

The complete enumeration of isotopy classes of Morse polynomials of degree $\leq 4$ 
with non-discriminant principal homogeneous parts in two variables is given in \cite{rigidplane}. We now move to the next dimension and consider degree three polynomials in ${\mathbb R}^3$. Our method is similar to  the one used in \cite{rigidplane}, but the three-dimensional situation has some considerable complications, see e.g. Remark \ref{spec}.

The parallel problem of classifying nonsingular surfaces of degree three was solved in \cite{schlaf}, \cite{klein}, \cite{zeuten} for surfaces in ${\mathbb R}P^3$, while a solution of its affine variant has only recently been published, see \cite{FK}. It turns out that the classification lists of isotopy classes of Morse polynomials of a certain degree and number of variables are more complicated than the lists of such polynomials with nonsingular zero sets. \smallskip

\begin{definition} \rm
\label{fdef}
A polynomial $({\mathbb C}^3, {\mathbb R}^3) \to ({\mathbb C}, {\mathbb R})$ is {\em Morse} if all its {\em real} critical points are Morse.
It is {\em strictly Morse} if, in addition, all its critical values at {\em real} critical points are different. The principal homogeneous part of a polynomial is {\em non-discriminant} if its zero set defines a smooth curve in ${\mathbb C}P^2$. A polynomial is {\em generic} if it has only Morse critical points in ${\mathbb C}^3$, all of their critical values are different from each other and from 0, and its principal homogeneous part defines a smooth curve in ${\mathbb C}P^2$. The polynomials lying in the same connected component of the space of Morse (or strictly Morse, or generic) polynomials of a given degree are called {\em isotopic} to each other in the corresponding space; such connected components are called {\em isotopy classes}.
\end{definition} 
\smallskip

The space of all degree three polynomials ${\mathbb R}^3 \to {\mathbb R}$ with 
non-discri\-mi\-nant principal homogeneous parts is divided into two connected subsets characterized by the number of components of the cubic curves in ${\mathbb R}P^2$ defined by these principal parts. We call the polynomials with one (respectively, two) such components {\em polynomials of type} $\Xi_1$ (respectively, $\Xi_2$).

Each of these two subsets is further divided into several connected components by the variety of non-Morse polynomials.

\begin{table}
\caption{Numbers of isotopy classes of Morse functions of types $\Xi_1$ and $\Xi_2$}
\label{tab0001}
\begin{tabular}{|l|c|c|c|c|c|c|}
\hline
{\rm Number of real critical points} & 0 & 2 & 4 & 6 & 8 & Total \\
\hline
{\rm Number of isotopy classes}, $\Xi_1$ & 1 & 3 & 4 & 5 & 8 & 21 \\
\hline
{\rm Number of isotopy classes,} $\Xi_2$ & 0 & 1 & 3 & 4 & 8 & 16 \\
\hline
\end{tabular}
\end{table}

\begin{theorem}
\label{mthm}
The numbers of isotopy classes of Morse polynomials of types $\Xi_1$ and $\Xi_2$ with different numbers of real critical points are given in Table \ref{tab0001}.
\end{theorem}

These isotopy classes of type $\Xi_1$ (respectively, $\Xi_2$) will be explicitly presented in \S \ref{p81} (respectively, \S \ref{p82}). 

\begin{theorem}[see \S \ref{strict}]
\label{mthm2}
There are exactly 842 $($respectively, 1416$)$
 isotopy classes of {\em strictly} Morse polynomials of type $\Xi_1$ $($respectively, $\Xi_2)$ with eight real critical points. Each of these isotopy classes is homotopy equivalent to the group $\mbox{SL}(3, {\mathbb R})$.
\end{theorem}

\subsection{Local theory of caustics}
\label{caus}
\
The principal parts of non-dis\-cri\-mi\-nant homogeneous polynomials of degree three in ${\mathbb R}^3$ can be considered as functions with isolated critical points at the origin. The space of all such polynomials with a given principal part $f_3$ is a deformation of the function singularity $f_3:$ the parameters of this deformation are the coefficients of the lower-degree monomials.   If our polynomial $f$ splits into homogeneous summands as $f_3 + f_2 + f_1+f_0$, then all functions from the one-parametric family $t^3f(t^{-1}x, t^{-1}y, t^{-1}z) \equiv f_3 + t f_2 + t^2 f_1 +t^3 f_0,$ $t >0$, have the same topological type. Therefore, studying the isotopy classes of arbitrary polynomials  of class $\Xi_1$ or $\Xi_2$  we can assume that these polynomials are arbitrarily small perturbations of their principal parts $f_3$. Thus, our study includes enumerating the local components of the complements of caustics of deformations of function  singularities of class $P_8$. 

The geometry of the complements of caustics of the simplest singularity classes $A_2$, $A_3$, $A_4$ and $D_4^{\pm}$ was studied by R.~Thom and V.~Arnold in the context of catastrophe theoretical problems occurring in biology and optics, see \cite{thom}, \cite{thom2}, \cite{Acongr}, \cite{Kluwer}. 
The connected components of these complements of all {\em simple} function singularities were listed in \cite{sed}, \cite{sede}, and \cite{jsing}.  
The next natural set of singularity classes, the {\em parabolic} ones, consists of the classes $P_8$, $X_9$, and $J_{10}$. The problem of the topological classification of small Morse  perturbations of these singularities is solved respectively in the present work, in \cite{rigidplane}, and in \cite{Vj10}.

\subsection{Normal forms of homogeneous cubic polynomials on ${\mathbb R}^3$}
The non-discriminant principal homogeneous part of a third-degree polynomial on ${\mathbb R}^3$
can always be reduced to the Hessian normal form
\begin{equation}
\label{symnf}
\frac{1}{3}\left(x^3 + y^3 + z^3\right) - A x y z , \qquad A \neq 1 ,
\end{equation}
by a linear transformation of ${\mathbb R}^3$,
see \cite{BM}. The parameter $A$ in (\ref{symnf}) is below (respectively, above) 1 for polynomials of type $\Xi_1$ (respectively, $\Xi_2$). We will also use the Newton--Weierstrass normal form 
\begin{equation}
\label{wei}
x^3 + a x z^2 + b z^3 - y^2 z 
\end{equation}
of these homogeneous parts,
where the polynomial $t^3 + a t + b$ has a single real root for polynomials of type $\Xi_1$ and three roots for $\Xi_2$. 

\subsection{Morse complex} 
\label{morcom}

Every polynomial of type $\Xi_1$ or $\Xi_2$ has at most eight real critical points. Moreover, the sum of the Milnor numbers of all its critical points in ${\mathbb C}^3$ is equal to eight. 

For any Morse polynomial of type $\Xi_1$, the Morse complex defined by its real critical points is acyclic. In other words, if $m$ (respectively, $M$) is a real number lower (respectively, higher) than all real critical values of such a polynomial $f$, then the relative homology group $$H_*(f^{-1}((-\infty,M]), f^{-1}((-\infty, m])) $$ is trivial in all dimensions. Indeed, these homology groups are the same for all $f$ from the connected set of polynomials of type $\Xi_1$; this set contains in particular the polynomials 
\begin{equation}
\label{s}
x^3 + y^3 + z^3 + \varepsilon (x + y + z),
\end{equation} which have no real critical points if $\varepsilon >0$.

For every Morse polynomial of type $\Xi_2$, the corresponding Morse complex has non-trivial homology groups $H_1 \simeq H_2 \simeq {\mathbb Z}$, in particular every such polynomial has at least one critical point of Morse index 1 and one of index 2.

In both cases, the numbers of critical points with even and odd Morse indices are always equal to each other.
\smallskip

The simplest invariant of the isotopy classes of Morse polynomials with non-discriminant principal homogeneous parts is the following.
 
\begin{definition} \rm
The {\it passport} of a Morse polynomial ${\mathbb R}^3 \to {\mathbb R}$ is the set of numbers of its critical points with Morse indices 
equal to 0, 1, 2 and 3.
\end{definition}
\smallskip

For example, the passport of the polynomial (\ref{s}) is equal to $(0, 0, 0, 0)$ if $\varepsilon >0 $ and to $(1, 3, 3, 1)$ if $\varepsilon <0$. 

This invariant is not sufficient to separate all isotopy classes of degree three Morse polynomials. In \S \ref{minv} we describe much stronger invariants that separate all these classes.

\begin{proposition}
\label{no2max}
A polynomial of degree 3 cannot have more than one Morse minimum point or more than one Morse maximum point.
\end{proposition}

\noindent
{\it Proof.} The restriction of such a polynomial to the line through such two points would be a degree $\leq 3$ polynomial in one variable with two Morse minima or two maxima. \hfill $\Box$

\subsection{Splitting function singularities of type $P_8$ into pairs of critical points}

A large part of our examples of isotopy classes of Morse polynomials are related to polynomials of type $\Xi_1$ or $\Xi_2$ having two real critical points with the sum of Milnor numbers equal to eight. These polynomials appear in the singularity theory of smooth functions as obstacles to the equisingularity of bifurcation sets along $\mu = \mbox{const}$ strata.

\begin{definition} \rm
A {\em function singularity} is a smooth function germ $\varphi: {\mathbb R}^3 \to {\mathbb R}$ at a point where its differential vanishes.
Two function germs, $\varphi$ at point $a$ and $\psi$ at point $b$, are {\it equivalent} if there exists a local diffeomorphism $T:({\mathbb R}^3, a) \to ({\mathbb R}^3,b)$ of the argument space such that $\varphi \equiv \psi \circ T$ in a neighborhood of the point $a$. 
\end{definition} 

 Let $f$ be a polynomial of type $\Xi_1$ or $\Xi_2$ with a pair of critical points with the sum of the Milnor numbers equal to eight. These critical points are necessarily {\em simple} in the terminology of \cite{AVGZ1}, since all function singularities with Milnor numbers less than eight are simple. 
\begin{table}
\caption{Real simple and $P_8$ singularities in three variables (all signs $\pm$ are independent)}
\label{t1}
\begin{tabular}{|l|l|l|}
\hline
Notation & Normal form & Restriction \\ 
\hline 
$A_{2k-1}$ & $\pm x^{2k} \pm y^2 \pm z^2$ & $k \ge 1$ \\ 
$A_{2k}$ & $x^{2k+1} \pm y^2 \pm z^2 $ & $k \ge 1$ \\
[3pt]
$D_{2k}^+$ & $x^2y \/+ y^{2k-1} \pm z^2$ & $k \ge 2$ \\ [4pt]
$D_{2k}^-$ & $x^2y \/- y^{2k-1} \pm z^2$ & $k \ge 2$ \\ [4pt]
$D_{2k-1}$ & $\pm (x^2y \/+ y^{2k-2})\pm z^2$ & $k \ge 3$ \\ [3pt]
$E_6$ & $x^3 \pm y^4 \pm z^2$ & \cr
$E_7$ & $x^3 + x y^3 \pm z^2$ & \cr
$E_8$ & $x^3 + y^5 \pm z^2$ & \cr
\hline
$P_8^1$ & $\frac{1}{3}(x^3+y^3+z^3) - A x y z$ & $A<1$ \\
$P_8^2$ & $\frac{1}{3}(x^3+y^3+z^3) - A x y z$ & $A>1$ \\
\hline
\end{tabular} 
\end{table}
(See Table \ref{t1} for the list of all real simple singularity classes of functions ${\mathbb R}^3 \to {\mathbb R}^1$ and the normal forms, to which the function germs of these classes can be reduced by choosing local coordinates and adding constant functions). The one-parameter family of polynomials \
$t^3 f(t^{-1}x, t^{-1}y, t^{-1}z), \ t \to +0, $ \ all having the same homogeneous part of degree three and the same collection of classes of critical points, then tends to a function with a singularity of class $P_8$ (i.e., with vanishing quadratic part of the Taylor decomposition) at the origin. 

The set $\{P_8\}$ of functions with singularities of class $P_8$ has the same codimension (equal to 6) in the space of all functions $M^3 \to {\mathbb R}^1,$ as the set of functions with two simple critical points with the sum of Milnor numbers equal to eight. Therefore, the functions of class $P_8$ that can be approached by the latter set form a subset of positive codimension in the entire set $\{P_8\}$. This property is invariant under the equivalences of function germs, so it is satisfied for $P_8$-type polynomials with only finitely many values of the module $A$ of the normal form (\ref{symnf}). The stratification of the space of functions according to the singularity types of their critical points behaves in an exceptional way at these polynomials of type $\{P_8\}$ compared to the generic polynomials of this type.

 We prove the following two theorems.

\begin{theorem}
\label{adj1}
A. There exist polynomials of type $\Xi_1$ with pairs of real critical points of the following classes: \ $E_6+A_2$, \ $D_5+A_3$, \ $A_4+A_4,$ \ $D_4^+ + D_4^+$, \ $D_4^- + D_4^-$.

B. There are no polynomials of type $\Xi_1$ with the pairs of real critical points of the following classes: \ $A_7 + A_1$, \ $D_7 + A_1$, \ $E_7 + A_1$, \ $A_6 + A_2$, \ $D_6^+ + A_2$, \ $D_6^- + A_2$, \
$A_5 + A_3$, \ $A_4 + D_4^+$, \ $A_4 + D_4^-$, \ $D_4^+ + D_4^-$. 
\end{theorem}

\begin{theorem}
\label{adj2}
A. There exist polynomials of type $\Xi_2$ with pairs of real critical points of classes $A_5 + A_3$ and $D_5 + A_3$.

B. There are no polynomials of type $\Xi_2$ with the pairs of real critical points of the following classes: \ $A_7 + A_1$, \ $D_7 + A_1$, \ $E_7 + A_1$, \ $A_6 + A_2$, \ $D_6^+ + A_2$, \ $D_6^- + A_2$, \ $E_6 + A_2$, \ $A_4 + A_4$, \ $A_4 + D_4^+$, \ $A_4 + D_4^-$, \ $D_4^+ + D_4^+$, \ $D_4^+ + D_4^-$, 
 \ $D_4^- + D_4^-$. 
\end{theorem}

\begin{remark} \rm
These two theorems are a ``real'' analog of the result of \cite{wall},  \cite{Jaw2} on the decomposition of the complex function singularities of type $P_8$ into pairs of critical points with the sum of Milnor numbers equal to eight. The real version of this result consists of a larger number of separate statements due to different real forms of these critical points.
\end{remark}

\subsection{On the structure of the work}

In \S~\ref{minv} we introduce the main invariants of isotopy classes of Morse polynomials of types $\Xi_1$ and $\Xi_2$ and prove their basic properties.

In \S\S~\ref{p81} and \ref{p82} we describe all values of these invariants, prove their completeness and realize them by polynomials. In \S~\ref{strict} we derive from these results a description of all isotopy classes of {\em strictly} Morse polynomials of types $\Xi_1$ and $\Xi_2$ with only real critical points, in particular find the numbers of them. In \S~\ref{pt4} we prove Theorems \ref{adj1} and \ref{adj2}.

\subsection*{A notation} The symbol $\Box$ indicates either the conclusion of a proof or the absence of a proof, as in the case of immediate corollaries or statements supplied with references to the works where they are proven. 

\section{Main invariants}
\label{minv}

\subsection{Virtual morsifications and set-valued invariants}

\unitlength 1.00mm
\linethickness{0.4pt}

\begin{figure}
\begin{center}
\begin{picture}(60,40)
\put(0,20){\line(1,0){60}}
\put(6.6,20){\circle*{1.5}}
\put(16.6,20){\circle*{1.5}}
\put(41.6,20){\circle*{1.5}}
\put(51.6,20){\circle*{1.5}}
\bezier{200}(30,20)(10,35)(5,20)
\bezier{100}(30,20)(18,25)(15,20)
\bezier{100}(30,20)(35,25)(40,20)
\bezier{200}(30,20)(40,35)(50,20)
\bezier{15}(5,20.2)(6,20.2)(6.5,20.2)
\bezier{15}(15,20.2)(16,20.2)(16.5,20.2)
\bezier{15}(40,20.2)(41,20.2)(41.5,20.2)
\bezier{15}(50,20.2)(51,20.2)(51.5,20.2)
\put(13,35){\circle*{1.5}}
\put(13,5){\circle*{1.5}}
\put(45,37){\circle*{1.5}}
\put(45,3){\circle*{1.5}}
\bezier{200}(30,20)(25,32)(13,35)
\bezier{200}(30,20)(25,8)(13,5)
\bezier{200}(30,20)(35,30)(45,37)
\bezier{200}(30,20)(35,10)(45,3)
\put(6,17){{\tiny 1}}
\put(16,17){{\tiny 2}}
\put(41,17){{\tiny 3}}
\put(51,17){{\tiny 4}}
\put(12,32){{\tiny 5}}
\put(44,34){{\tiny 6}}
\put(12,7){{\tiny 7}}
\put(44,5){{\tiny 8}}
\end{picture} \qquad 
\begin{picture}(60,40)
\put(0,20){\line(1,0){60}}
\put(6.6,20){\circle*{1.5}}
\put(16.6,20){\circle*{1.5}}
\put(41.6,20){\circle*{1.5}}
\put(51.6,20){\circle*{1.5}}
\bezier{200}(30,20)(10,35)(5,20)
\bezier{100}(30,20)(18,25)(15,20)
\bezier{100}(30,20)(35,25)(40,20)
\bezier{200}(30,20)(40,35)(50,20)
\bezier{15}(5,20.2)(6,20.2)(6.5,20.2)
\bezier{15}(15,20.2)(16,20.2)(16.5,20.2)
\bezier{15}(40,20.2)(41,20.2)(41.5,20.2)
\bezier{15}(50,20.2)(51,20.2)(51.5,20.2)
\put(13,35){\circle*{1.5}}
\put(13,5){\circle*{1.5}}
\put(45,37){\circle*{1.5}}
\put(45,3){\circle*{1.5}}
\bezier{200}(30,20)(25,32)(13,35)
\bezier{200}(30,20)(25,8)(13,5)
\bezier{200}(30,20)(22,30)(10,31)
\bezier{70}(10,31)(5,31.5)(5,34)
\bezier{70}(5,34)(5,37)(10,37)
\bezier{200}(10,37)(35,37)(45,37)
\bezier{200}(30,20)(22,10)(10,9)
\bezier{100}(10,9)(5,8.5)(5,6)
\bezier{100}(5,6)(5,3)(10,3)
\bezier{200}(10,3)(35,3)(45,3)
\put(6,17){{\tiny 1}}
\put(16,17){{\tiny 2}}
\put(41,17){{\tiny 3}}
\put(51,17){{\tiny 4}}
\put(12,32.3){{\tiny 6}}
\put(44,34){{\tiny 5}}
\put(12,6.4){{\tiny 8}}
\put(44,5){{\tiny 7}}
\end{picture}
\end{center}
\caption{Standard systems of paths}
\label{standd}
\end{figure}

If $f$ is a generic polynomial of degree three (see Definition \ref{fdef}), then the set $V_f \equiv f^{-1}(0) \subset {\mathbb C}^3$ is a smooth complex surface homotopy equivalent to the wedge of eight two-dimensional spheres (see \cite{M}). The homology group $H_2(V_f) \simeq {\mathbb Z}^8$ is generated by {\em vanishing cycles} (see e.g. \cite{AVGZ2}) defined by some non-intersecting paths in ${\mathbb C}^1$, connecting the non-critical value 0 with all critical values of $f$, see Fig.~\ref{standd}. We choose these paths and the orientations of the vanishing cycles in a standard way, see \S~V.1.6 of \cite{RISL}. 

Namely, the canonical orientation of cycles vanishing at real critical points is induced by a fixed orientation of ${\mathbb R}^3$ as follows. Let $\xi_1, \xi_2, \xi_3$ be real local coordinates at such a point, in which the function is equal to \begin{equation}
\label{orie}
\xi_1^2+ \dots +\xi_k^2-\xi_{k+1}^2 - \dots - \xi_3^2 + c_i, 
\end{equation} where $k \in \{0, 1, 2, 3\}$ and $c_i$ is the critical value; we choose these coordinates so that the orientation of ${\mathbb R}^3$ defined by the frame $\left(\frac{\partial}{\partial \xi_1}, \frac{\partial}{\partial \xi_2}, \frac{\partial}{\partial \xi_3}\right)$ coincides with the fixed one. The vanishing cycle $\Delta_i (c_i-\varepsilon) \subset f^{-1}(c_i - \varepsilon)$ for a very small $\varepsilon >0$ can be realized by a sphere of radius $\sqrt{\varepsilon}$ in the real 3-subspace of the corresponding chart of ${\mathbb C}^3$, which is spanned over ${\mathbb R}$ by the vectors 
\begin{equation}
\label{tsp}
i \frac{\partial}{\partial \xi_1}, \dots, i \frac{\partial}{\partial \xi_k}, \frac{\partial}{\partial \xi_{k+1}},\dots, \frac{\partial}{\partial \xi_3}
\end{equation}
 and is oriented by this sequence of vectors. The canonical orientation of this cycle is defined by an arbitrary tangent 2-frame such that the 3-frame $\{\mbox{grad } f, \mbox{this 2-frame}\}$ defines this orientation of the space (\ref{tsp}). The cycle $\Delta_i \subset f^{-1}(0)$ is obtained from the oriented cycle $\Delta_i(c_i-\varepsilon)$ by the covering homotopy transportation over the path connecting 0 and $c_i - \varepsilon$ in the upper half-plane in the area where the imaginary part is smaller than the absolute values of the imaginary parts of all non-real critical values of our polynomial. 

Any two paths leading to complex conjugate non-real critical values should be symmetric to each other with respect to the real axis. The orientation of cycles vanishing at non-real critical points is determined by the condition that the complex conjugation in ${\mathbb C}^3$ takes such cycles defined by conjugate paths to each other, and by some lexicographic maximization of intersection indices of these vanishing cycles over all systems of orientations satisfying the previous conditions. 

We also canonically order the critical values and the corresponding vanishing cycles: first the real critical values in the ascending order, then the critical values in the upper half-plane in the order of the angles of the corresponding paths with the axis of negative numbers, and finally the critical values with negative imaginary parts in the order repeating the order of their complex conjugates.

\begin{definition} \rm
\label{vmf}
A {\em virtual Morse function} associated with a generic Morse polynomial $f:({\mathbb C}^3, {\mathbb R}^3) \to ({\mathbb C}, {\mathbb R})$ of degree three is a collection of its topological data consisting of

a) the $8 \times 8$ matrix of the intersection indices in $V_f$
of canonically ordered and oriented vanishing cycles $\Delta_i \in H_2(V_f)$ corresponding to all critical values of $f$,

b) the string of intersection indices in $V_f$ of these vanishing cycles with the naturally oriented set $V_f \cap {\mathbb R}^3$ of real points, 

c) the string of {\em parities} of Morse indices of all real critical points of $f$, and

d) the number of negative and positive critical values of $f$.
\end{definition}

\noindent
{\bf Example.} A virtual Morse function associated with a polynomial of type $\Xi_1$ with six real critical points is given in (\ref{virtuM}) (left). The vertical lines in this notation indicate the last element of the virtual function: the numbers of negative and positive critical values are equal to 4 and 2, respectively.
\begin{equation}
\label{virtuM}
\begin{array}{|cccc|cc|cc|}
\hline
 $-2$ & 0 & 0 & 1 & 0 & 0 & 0 & $-1$ \\
 0 & $-2$ & 0 & 1 & 0 & 1 & 1 & 0 \\
 0 & 0 & $-2$ & 1 & 1 & 0 & 1 & 0 \\
 1 & 1 & 1 & $-2$ & 0 & 0 & $-1$ & 1 \\
 0 & 0 & 1 & 0 & $-2$ & 0 & 0 & 0 \\
 0 & 1 & 0 & 0 & 0 & $-2$ & 0 & 0 \\
 0 & 1 & 1 & $-1$ & 0 & 0 & $-2$ & $-1$ \\
 $-1$ & 0 & 0 & 1 & 0 & 0 & $-1$ & $-2$ \\
\hline
 1 & 1 & 1 & 0 & $-2$ & $-2$ & $-1$ & $-1$ \\
\hline
 o & o & o & e & e & e & & \\
\hline
\end{array} \qquad \qquad 
\begin{array}{cccccccc}
& & & & & & 6 & 4 \\
 $-2$ & 0 & 0 & 1 & 0 & 0 & 0 & $-1$ \\
 0 & $-2$ & 0 & 1 & 0 & 1 & 1 & 0 \\
 0 & 0 & $-2$ & 1 & 1 & 0 & 1 & 0 \\
 1 & 1 & 1 & $-2$ & 0 & 0 & $-1$ & 1 \\
 0 & 0 & 1 & 0 & $-2$ & 0 & 0 & 0 \\
 0 & 1 & 0 & 0 & 0 & $-2$ & 0 & 0 \\
 0 & 1 & 1 & $-1$ & 0 & 0 & $-2$ & $-1$ \\
 $-1$ & 0 & 0 & 1 & 0 & 0 & $-1$ & $-2$ \\
 & & & & & & & \\
 1 & 1 & 1 & $-2$ & 0 & 0 & 1 & $-1$ \\
 1 & 1 & 1 & 0 & $-2$ & $-2$ & $-1$ & $-1$ \\
& & & & & & & \\
 $1$ & $1$ & $1$ & $-1$ & $-1$ & $-1$ & 0 & 0 \\
\end{array} 
\end{equation}

\begin{remark} \rm
If a generic polynomial has more than one pair of non-real critical values, then there is more than one virtual Morse function associated with it, because in this case there is no canonical choice of a system of paths: see two parts of Fig.~\ref{standd}.
\end{remark}

\begin{definition} \rm 
A {\em critical point of a virtual Morse function} is any column of its data set as in (\ref{virtuM}) (left), i.e., the corresponding column of the intersection matrix, the intersection index of the corresponding vanishing cycle with the real space, and the parity of the Morse index or the information that the point is not real.
\end{definition}

\begin{definition} \rm
\label{elsur}
{\em Elementary virtual surgeries} of virtual Morse functions involve six transformations of their data, modeling the basic topological surgeries of the corresponding real generic polynomials, namely
\begin{itemize}
\item[s1, s2] \ \ collision of two neighboring real critical values, after which the corresponding two critical points either (s1) meet and leave the real space, or (s2) change the order in ${\mathbb R}^1$ of their critical values, 

\item[s3, s4] \ \ collision of two complex conjugate critical values at a point of the line ${\mathbb R}^1$, after which the corresponding critical points of $f_\lambda$ either (s3) meet at a real point and enter the real space, or (s4) miss each other while the imaginary parts of their critical values change their signs, 

\item[s5, s6] \ \ jumps of real critical values up (s5) or down (s6) through 0, \\ \hspace*{4cm} and additionally

\item[s7] \ \ specifically virtual transformations within the sets of virtual Morse functions associated with the same real Morse polynomials, caused by changes of the standard systems of paths going from 0 to non-real critical values (see two parts of Fig.~\ref{standd} and \S~2.6 in \cite{AVGZ2}).
\end{itemize}
\end{definition}

The results of all these virtual surgeries are determined by the data of the original virtual Morse functions. For a detailed description and explicit formulas of these standard flips of data, see \S~V.8--9 of \cite{RISL}. In particular, an attempt to perform the surgery $s1$ or $s2$ over the neighboring real critical values $v_i$ and $v_{i+1}$ begins by examining the intersection index $\langle \Delta_i, \Delta_{i+1}\rangle$ of the corresponding vanishing cycles. If this index is equal to $0$ then the surgery $s2$ occurs, if it is equal to 1 then the surgery $s1$ occurs, in all other cases the surgery fails. Similarly, the collision of two complex conjugate critical values at a point of ${\mathbb R}^1 \setminus \{0\}$ not separated from 0 by other critical values will follow scenario $s4$ if the intersection index of the corresponding vanishing cycles is equal to 0, scenario $s3$ if this index is equal to $1$ or $-1$, and will fail in all other cases; in the second case the sign of this intersection index allows us to predict the parities of the Morse indices of the newborn real critical points.

We will denote by $s1, \dots, s6$ both the real surgeries of real Morse functions and the corresponding elementary virtual surgeries.

\begin{definition} \rm
An (abstract) {\em virtual Morse function} of type $\Xi_1$ or $\Xi_2$ is any collection of data as in Definition \ref{vmf} (i.e., a matrix, two strings, and two numbers) obtained from a virtual Morse function associated with any real polynomial of that type by arbitrary chains of elementary virtual surgeries. 

The {\em formal graph} of type $\Xi_1$ or $\Xi_2$ is the graph whose vertices correspond to all virtual Morse functions of this type, and two vertices are connected by an edge if and only if the corresponding virtual Morse functions are obtained from each other by an elementary virtual surgery.

A {\em virtual component} of the set of all virtual Morse functions of type $\Xi_1$ or $\Xi_2$ is any of its maximal subsets, all elements of which can be obtained from each other by chains of elementary virtual surgeries not including the virtual models of the collisions of critical points, i.e. the surgeries s1 and s3 from Definition \ref{elsur} (but possibly including the surgeries s2 and s4, i.e., the collisions of critical values not related to the collisions of the corresponding critical points).
\end{definition}

\begin{remark} \rm
Two vertices of the formal graph can be connected by several edges, each corresponding to a different surgery. Also, this graph may contain edges with coinciding endpoints.
\end{remark}

Our main method consists of a combinatorial study of formal graphs and interpreting the results in the terms of real functions associated with virtual ones. 

\begin{remark} \rm
\label{spec}
The analogous technique for functions in two variables, used in particular in \cite{jsing},  \cite{rigidplane}, and \cite{Vj10}, is considerably easier, because in that case we can predict not only the parities of the Morse indices of newborn real critical points after surgery $s3$, but the integer indices themselves. Consequently, we now use a different set of programs for the combinatorial study of virtual Morse functions, and make additional efforts to reconstruct the integer Morse indices from the combinatorial data, see e.g. Proposition \ref{lemfun}.
\end{remark}

\begin{theorem}
\label{propmain}
1. If two generic polynomials of the same type $\Xi_1$ or $\Xi_2$ are isotopic to each other in the space of Morse polynomials, then all virtual Morse functions associated with these polynomials belong to the same virtual component. 

2. For any generic Morse polynomial $f$ of type $\Xi_1$ or $\Xi_2$, a virtual Morse function $\varphi$ associated with it, and a virtual Morse function $\tilde \varphi$ connected with $\varphi$ by an edge of the formal graph, there exists a generic Morse polynomial $\tilde f$ and a path in the space $\Xi_1$ or $\Xi_2$ connecting $f$ and $\tilde f$ and containing only one non-generic point at which it experiences a
standard surgery  of the same type as the edge $[\varphi, \tilde \varphi]$.
\end{theorem}

\begin{corollary}
\label{cormain}
1. Any virtual Morse function of the type $\Xi_1$ or $\Xi_2$ is associated with a real polynomial of that type. 

2. If two virtual Morse functions belong to the same virtual component, then they are associated with two generic Morse polynomials which are isotopic to each other in the set of Morse polynomials. \hfill $\Box$
\end{corollary}

\noindent
{\it Proof of Theorem \ref{propmain}.} Statement 1 is tautological: we connect our two polynomials by a smooth path in the space of Morse polynomials with non-discriminant principal parts and apply to the associated virtual Morse functions all surgeries occurring at the intersections of this path with the variety of non-generic polynomials.

If the edge $[\varphi, \tilde \varphi]$ in statement 2 encodes a 
surgery of type $s5$ or $s6$, then its realization in the space $\Xi_1$ or $\Xi_2$ is achieved by adding constants to $f$, and the surgery $s7$ is completely virtual and does not require any modification of $f$ at all. Therefore, we only need to realize all elementary surgeries of the types $s1$---$s4$ which are related to collisions of critical values. 

By a choice of affine coordinates in ${\mathbb R}^3$, we can assume that $f$ belongs to the standard miniversal deformation of functions of class $P_8$, i.e. to the family of polynomials
\begin{equation}
f_{A, \lambda} \equiv\frac{1}{3}(x^3 + y^3 + z^3) - A x y z + \lambda_1 + \lambda_2 x + \lambda_3 y + \lambda_4 z + \lambda_5 x y + \lambda_6 x z + \lambda_7 y z
\label{vers0}
\end{equation}
(see e.g. \cite{AVGZ1}, \cite{Jaw2}),
where $A$ and $\lambda_1, \dots , \lambda_7$ are {\em complex} parameters of the deformation (taking real values for $f$), $A$ can take any values except for only three, $1$ and $\pm e^{2\pi i/3},$ that correspond to discriminant homogeneous polynomials. Indeed, the normal form (\ref{symnf}) of the principal homogeneous part can be achieved according to \cite{BM}, and the monomials $x^2, y^2, z^2$ in $f$ can be eliminated by choosing the origin in ${\mathbb R}^3$. Let $\Theta$ denote the parameter space $\left({\mathbb C}^1 \setminus \{1, e^{\pm 2\pi i/3}\} \right) \times {\mathbb C}^7$ of this family.

Our construction is based on the results of \cite{Jaw}, \cite{Jaw2} and exploits the {\it Lyashko--Looijenga map} $\Lambda: \Theta \to \mbox{Sym}^8({\mathbb C}^1) \equiv {\mathbb C}^8 / S(8)$ (see \cite{Lo0}, \cite{Lo1}), which sends any polynomial (\ref{vers0}) with non-dis\-cri\-mi\-nant principal part to the unordered collection of its critical values (taken with multiplicities in the case of not strictly Morse polynomials). 

\begin{lemma}[see \cite{Jaw}]
This map $\Lambda$ is a covering over the smooth space $B({\mathbb C}^1, 8) \subset \mbox{\rm Sym}^8({\mathbb C}^1) $ of all unordered collections of eight {\em pairwise distinct} points in ${\mathbb C}^1$, i.e., in the restriction to the space of polynomials $($\ref{vers0}$)$ with eight different critical values. \hfill $\Box$
\end{lemma}

\begin{corollary}
\label{reallem}
If $f \in \Theta$ is a generic real polynomial, then the subspace $\Theta \cap {\mathbb R}^8 \subset \Theta$ of real polynomials of the form $($\ref{vers0}$)$ coincides in its neighborhood with the set of polynomials of class $\Theta$ whose sets of critical values are invariant under the complex conjugation. \hfill $\Box$
\end{corollary}

Let $I:[0,1] \to \mbox{Sym}^8({\mathbb C}^1)$ be an arbitrary piecewise algebraic path such that 
\begin{itemize}
\item
$I(0)$ is the collection of critical points of $f$, 
\item
$I(\tau) \in B({\mathbb C}^1, 8)$ for $\tau \in [0,1)$, 
\item some six points of the configuration $I(0)$ are fixed for all $\tau \in [0,1]$, and two remaining points of $I(\tau)$ collide as in the prescribed surgery when $\tau$ tends to $1$; 
\item
if we realize a surgery of type $s3$ or $s4$ then these two points should be mutually complex conjugate along the path, and in the case of surgeries $s1$ and $s2$ they should be real all the time. 
\end{itemize}

Let us try to lift the path $I$ to a map $\tilde I: [0,1] \to \Theta$ such that $\tilde I (0) = f$ and $I \equiv \Lambda \circ \tilde I$.
By the covering homotopy property of the covering $\Lambda$ over $B({\mathbb C}^1,8)$, this map (if it exists) is unique. Due to the algebraicity of the path $I$ and the map $\Lambda$, the limit 
\begin{equation}\lim_{\tau \to 1}\tilde I(\tau)
\label{limlim}\end{equation}
 in the compactification ${\mathbb C}P^1 \times {\mathbb C}P^7$ of the space $\Theta$ is well defined. By Proposition 2 of \cite{Jaw2}, this limit point belongs to the space $\left({\mathbb C}^1 \setminus \{1, e^{\pm 2\pi i/3}\}\right) \times {\mathbb C}P^7$. Indeed, the limit partition of the Dynkin diagram is always elliptic: if we perform one of the surgeries s2 or s4 then it splits into eight single vertices, in the case of surgery s3 it has only one non-trivial subdiagram 
\begin{picture}(15,3)
\put(3,1){\circle*{1}}
\put(12,1){\circle*{1}}
\put(3.5,1){\line(1,0){8}}
\end{picture} or 
\begin{picture}(15,3)
\put(3,1){\circle*{1}}
\put(12,1){\circle*{1}}
\bezier{30}(4.5,1)(5.5,1)(6.5,1)
\bezier{30}(8.5,1)(9.5,1)(10.5,1)
\end{picture}, 
and in the case of surgery s1 only the first of these subdiagrams.
 Therefore, the $A$-coordinates of the points $\tilde I(\tau) \in \Theta$ for $\tau<1$ sufficiently close to $1$ lie in some compact disk $\bar U \subset {\mathbb C}^1 \setminus \{1, e^{\pm 2\pi i/3}\}$. It follows from the connectedness of the Dynkin diagram (see \cite{Gab}) that there are no functions $f_{A, \lambda}$ with $A \in \bar U$ and $\lambda \neq 0$ such that all critical values of $f_{A, \lambda}$ are equal to 0. Due to the compactness of $\bar U$, this implies that there exists a number $\delta>0$ such that $\| \lambda \| \leq 1$ for all functions $f_{A, \lambda}$, $A \in \bar U$, whose critical values all lie in the $\delta$-neighborhood of 0. Therefore, the $\lambda$-coordinates of points $\tilde I(\tau)$ cannot tend to infinity when $\tau$ tends to $1$, and the limit point (\ref{limlim}) belongs to $\Theta$. This point is a polynomial with seven different critical values: it has a critical point of class $A_2$ if our surgery is of type $s1$ or $s3$, or two Morse critical points with the same real critical value in the case of surgery $s2$ or $s4$. By Corollary \ref{reallem}, all the polynomials $\tilde I(t),$ $t \in [0,1),$ are real, so the same is true for $\tilde I(1)$, and also this limit double critical value is real. The desired Morse perturbation of the corresponding real critical point of class $A_2$ (in the case of surgery $s1$ or $s3$) or the shift of the critical values of two complex conjugate Morse critical points (in the cases $s2$ and $s4$) can be performed within the real part of the family (\ref{vers0}) by the versality of this family.
 \hfill $\Box$

\subsection{Reconstructing integer Morse indices from virtual Morse functions}
Although the data of virtual Morse functions contain only the parities of the Morse indices of the critical points of the associated Morse polynomials, in many cases it is possible to reconstruct these indices themselves. 

\begin{proposition} 
\label{lemfun}
Suppose that the intersection index $\langle \Delta_i, \Delta_j \rangle$, $i<j$, of some two vanishing cycles corresponding to real critical points
of a generic polynomial $f$ of type $\Xi_1$ or $\Xi_2$ is equal to 1, the critical values of these critical points have the same sign, and there are no increasing chains of numbers $i=i_1 < i_2< \dots < i_k =j$, $k\geq 3$, such that all the intersection indices $\langle \Delta_{i_p}, \Delta_{i_{p+1}} \rangle$, $p=1, \dots, k-1$, are not equal to 0. Then the Morse index of the $i$th critical point of $f$ is by 1 less than that of the $j$th critical point.
\end{proposition}

\noindent
{\it Proof.} For any intermediate natural number $ l \in (i, j)$ either there are no increasing chains of numbers $i_p$ as above connecting $l$ and $j$ and such that $\langle \Delta_{i_p}, \Delta_{i_p+1} \rangle \neq 0$ for any pair of vanishing cycles which are neighbors in this sequence, or there are no such chains growing from $i$ to $l$. By continuously shifting up all critical values that satisfy the first property, and shifting down all critical values that satisfy the second property (but not both), we obtain a point configuration in ${\mathbb R}^1$ that contains no points between our two critical values. During this shift only the collisions of critical values corresponding to vanishing cycles with zero intersection indices occur, so the lifting of this shift to $\Theta$ starting from the polynomial $f$ as in the proof of Theorem \ref{propmain} can be performed within the set of Morse polynomials. Two critical points of the resulting Morse polynomial $\tilde f$, corresponding to our two critical values, have the same Morse indices as at the beginning, the intersection index of the corresponding vanishing cycles is still equal to 1, and there are
no other critical values between them. Consider another path in the space $\mbox{Sym}^8({\mathbb C}^1)$, along which all other critical values remain fixed and these two collide at an intermediate point. Lifting this path to $\Theta$ with the starting point $\tilde f$ tends to a polynomial with a critical point of class $A_2$, from which our two critical points are obtained by the standard Morse decomposition. \hfill $\Box$ \medskip

\noindent
{\bf Example.} \label{exaa} The Morse indices of the first three critical points of any polynomial associated with the virtual Morse function (\ref{virtuM}) are lower by 1 than the Morse index of the fourth point, so they are not equal to 3. Since they are odd, they can only be equal to 1.

\begin{theorem}
\label{refprop}
If two generic polynomials $f$, $\tilde f $ are associated with the same virtual Morse function, then either they are isotopic to each other in the space of generic polynomials, or their isotopy classes in this space 
are mapped into each other by any affine involution of the argument space ${\mathbb R}^3$ changing its orientation. \end{theorem}

Proofs of this theorem for polynomials of type $\Xi_1$ and $\Xi_2$ are given separately in subsections \ref{prorefprop} and \ref{prorefprop2}. 

\begin{definition} \rm
\label{defsv}
The {\em set-valued invariant} of a generic real Morse polynomial ${\mathbb R}^3 \to {\mathbb R}$ of degree three is the entire virtual component of any virtual Morse function associated with that polynomial. 
The invariant Card of such a polynomial is the number of elements of this virtual component.
\end{definition}
\smallskip

By Theorems \ref{propmain} and \ref{refprop}, the set-valued invariant is indeed an invariant of isotopy classes of Morse polynomials, and 
every possible value of this invariant (i.e., every virtual component of type $\Xi_1$ or $\Xi_2$) can be only taken on one or two such isotopy classes.

\begin{definition} \rm
\label{her}
A generic polynomial $f$ is {\em achiral} (respectively, {\em chiral}) if every affine involution of the argument space ${\mathbb R}^3$, which changes the orientation, takes its isotopy class $\{f\}$ of Morse polynomials to itself (respectively, to a different isotopy class). 
A virtual component is achiral (respectively, chiral) if it is the value of the set-valued invariant of achiral (respectively, chiral) Morse polynomials.
\end{definition}

\begin{remark} \rm
An analog of the set-valued invariant can be defined for isotopy classes of higher degree Morse polynomials, but in these cases there are no direct analogs of Theorems \ref{propmain} (except for its first statement) and \ref{refprop}. 
\end{remark} 

\begin{proposition}
\label{procount1}
There are exactly 6503 $($respectively, 9174$)$ virtual Morse functions associated with generic real polynomials ${\mathbb R}^3 \to {\mathbb R}$ of type $\Xi_1$ $($respectively, $\Xi_2)$. 
\begin{table}
\caption{Numbers of virtual Morse functions of types $\Xi_1$ and $\Xi_2$}
\label{tab01}
\begin{tabular}{|l|c|c|c|c|c|c|}
\hline
{\rm Number of real critical points} & 0 & 2 & 4 & 6 & 8 & Total \\
\hline
{\rm Number of virtual functions}, $\Xi_1$ & 297 & 390 & 515 & 1512 & 3789 & 6503 \\
\hline
{\rm Number of virtual functions,} $\Xi_2$ & 0 & 255 & 650 & 1897 & 6372 & 9174 \\
\hline
\end{tabular}
\end{table}
The numbers of virtual Morse functions of these two types
with different numbers of real critical points are given respectively in the second and the third rows of Table~\ref{tab01}. 
\end{proposition}

\noindent
{\it Proof}. These numbers were found using the computer program 

\noindent
{\footnotesize
\begin{verbatim}
https://drive.google.com/file/d/18QWyCZh6DVLEeu1nzDxMDwux7_InfHdV/view?usp=sharing
\end{verbatim}
\label{pro1}}

\noindent
described in \cite{RISL}, \S\S~V.8--9, and \cite{AGLV2}, and also used in \cite{Petr9}.
Namely, for any type $\Xi_1$ or $\Xi_2$ we start with the intersection matrix of the basic vanishing cycles (defined by a system of paths as in Fig.~\ref{standd}) of a certain Morse polynomial of this type having only real critical points. In the case $\Xi_1$ we take for this the polynomial $x^3 + y^3 +z^3$, whose intersection matrix can be computed by the method of A.~Gabrielov, see e.g. \S~II.2.3 in \cite{AGLV1}, and is represented by the Coxeter--Dynkin diagram shown in Fig.~\ref{72} below. For $\Xi_2$ the analogous matrix is computed in \S~7 of \cite{Petr9}. Then we insert all its non-zero numbers $\langle \Delta_i, \Delta_j \rangle$, $i<j$, as the values $C(i,j)$ at the end of the program, and also set there $\mbox{INDC}(i)=-1$ for all numbers $i$ such that the $i$th critical value has an even Morse index. Finally, the number of negative critical values should be indicated as the value of the parameter NPOZC in line 39 of the program.

From these initial data, the program first computes the intersection indices of the vanishing cycles with the set of real points, and thus obtains all the elements of the virtual Morse function associated with our polynomial.
Then it performs arbitrary chains of standard virtual surgeries and finds all virtual Morse functions in the corresponding formal graph, in particular it counts the numbers of these virtual Morse functions with all possible numbers of real critical points. \hfill $\Box$ \smallskip

\begin{corollary}
The sets of virtual Morse functions of types $\Xi_1$ and $\Xi_2$ have no common elements.
\end{corollary}

\noindent
{\it Proof.} Each virtual Morse function determines the entire set of all virtual Morse functions of the corresponding type, but these sets of types $\Xi_1$ and $\Xi_2$ are different. \hfill $\Box$ \smallskip

\begin{remark} \rm
In our program printouts, the virtual Morse functions are encoded in a
slightly different way than in this article: for example, the virtual Morse function given in (\ref{virtuM}) (left) appears as (\ref{virtuM}) (right). The letters ``o'' , ``e'' and the spaces in the bottom row on the left are encoded by the numbers $1$, $-1$ and 0 respectively, the numbers of {\it all real} and negative critical values can be read above the matrix on the right. 
It also contains an additional string (see the third row from the bottom), which is important in the theory of hyperbolic PDEs but can be ignored in our current problems. This string can be easily derived from the last two: for any column corresponding to a real critical value the sum of the numbers in the rows 2 and 3 from the bottom is equal to $2$ (respectively, to $-2$) if we have ``o'' (respectively, ``e'') in the last row; for vanishing cycles corresponding to non-real critical values with negative (respectively, positive) imaginary parts, the numbers in the rows 2 and 3 are equal (respectively, opposite). 

In PDE theory, the second and the third rows from bottom are called the {\em even} and {\em odd local Petrovskii cycles}, see \cite{RISL}. 

\end{remark}
\smallskip

 In \S\S~\ref{p81}, \ref{p82} below we describe all the virtual components of types $\Xi_1$ and $\Xi_2$, realize them by real polynomials associated with some their elements, and study their chirality. \smallskip

\begin{remark} \rm
It is not surprising that the number of virtual Morse functions with exactly $k$ real critical points is always divisible by $k+1$: each such virtual function is accompanied by $k$ others, which are obtained from it by iterated surgeries s5 and s6. Moreover, the same is true for all values of the Card invariant.
\end{remark}

\subsection{Up-down involution}
\label{invo}

The involution 
\begin{equation}
\label{invof}
f(x, y, z) \leftrightarrow -f(-x, -y, -z)
\end{equation}
acts on both spaces $\Xi_1$ and $\Xi_2$, preserving the sets of Morse polynomials and the principal homogeneous parts of these polynomials. 

This involution can be uniquely extended to virtual Morse functions: the set of virtual Morse functions associated with the polynomial $-f(-x, -y, -z)$ (or, equivalently, just with $-f$) is determined by the analogous set associated with $f$. In particular, there are the following relations between the vanishing cycles of such polynomials and between their intersection indices.

\begin{proposition}
\label{proud}
Suppose that a generic polynomial $f$ has $r$ real critical points and $s =(8-r)/2$ pairs of non-real ones, and the polynomial $\check f$ is obtained from $f$ by the involution $($\ref{invof}$)$. For any $i=1, \dots, r$, 
the vanishing cycle $\Delta_i$ of the initial polynomial $f$ and 
the vanishing cycle $\tilde \Delta_{r+1-i}$ of the polynomial $\check f$ are then obtained from each other by the composition of a$)$ central symmetry in ${\mathbb C}^3$ defining in particular a diffeomorphism $f^{-1}(0) \leftrightarrow \check f^{-1}(0)$, b$)$ complex conjugation, and c$)$ change of the orientation. 
For any $i=1, \dots, 2s,$ the central symmetry in ${\mathbb C}^3$ moves the vanishing cycle $\Delta_{r+i}$ of $f$ to the vanishing cycle $\pm \tilde \Delta_{9-i}$. The sign $\pm $ depends on $i$, but  for any $i \in \{1, \dots, s\}$ it is the same for $\Delta_i$ and $\Delta_{i+s}$.
\end{proposition}\smallskip

\begin{figure}
\unitlength 1mm
\begin{picture}(60,20)
\put(0,10){\line(1,0){60}}
\put(20,10){\circle*{1.3}}
\put(5,10){\circle*{1.3}}
\put(35,20){\circle*{1.3}}
\put(35,0){\circle*{1.3}}
\bezier{150}(30,10)(31,17)(35,20)
\bezier{150}(30,10)(31,3)(35,0)
\bezier{200}(30,10)(20,15)(18,10)
\bezier{300}(30,10)(10,22)(3,10)
\bezier{300}(30,10)(11,20)(6.5,10)
\put(28,6){\small 0}
\put(3.5,6){\small $c_i$}
 \end{picture} \qquad
\begin{picture}(60,20)
\put(0,10){\line(1,0){60}}
\put(40,10){\circle*{1.3}}
\put(55,10){\circle*{1.3}}
\put(25,20){\circle*{1.3}}
\put(25,0){\circle*{1.3}}
\bezier{150}(30,10)(29,17)(25,20)
\bezier{150}(30,10)(29,3)(25,0)
\bezier{200}(30,10)(37,7)(38.5,10)
\bezier{300}(30,10)(49,0)(53.5,10)
\bezier{300}(30,10)(49,20)(53.5,10)
\put(27.5,6.5){\small 0}
\put(52.6,6){\small $-c_i$}
 \end{picture} 
\caption{Up-down involution}
\label{udi}
\end{figure}

\noindent
{\it Proof} (see Fig.~\ref{udi}). Let $c_i$, $ i \in \{1, \dots, r\}$, be a real critical value of $f$. Together with the vanishing cycle $\Delta_i(c_i -\varepsilon) \subset f^{-1}(c_i - \varepsilon)$, the cycle $\Delta^+_i (c_i +\varepsilon) \subset f^{-1}(c_i +\varepsilon)$ can be defined in the same way. It is easy to compute (see Lemma 1 in \S V.3.1 of \cite{RISL}) that the transport over the arc $\{c_i + \varepsilon e^{it},\}, t \in [0,\pi],$ moves one of these two cycles to the other multiplied by $(-1)^{k}$, where $k$ is the positive inertia index of the second differential of the $i$-th critical point, see (\ref{orie}). By the construction, the central symmetry in ${\mathbb C}^3$ geometrically moves the cycle in $f^{-1}(0)$, which is transported from the cycle $\Delta^+_i(c_i+\varepsilon)$ over a path in the upper half-plane, into a cycle in $\check f^{-1}(0)$, which is transported from the cycle $\tilde \Delta_{r+1-i}(-c_i -\varepsilon) \subset \check f^{-1}(-c_i - \varepsilon)$ over a path in the lower half-plane. Moreover, this move preserves the canonical orientations of these cycles (in fact, one of them is multiplied by $(-1)^3 \times (-1) \equiv 1,$ where the first factor means the change of the orientation of ${\mathbb R}^3$ at the central symmetry, and the second factor results from a comparison of the directions of the gradients of our functions). The complex conjugation moves this cycle in $\check f^{-1}(0)$ to the cycle transported over a path in the upper half-plane from the cycle complex conjugate to $\tilde \Delta_{r+1-i}(-c_i-\varepsilon)$. The latter cycle is equal to the cycle $\tilde \Delta_{r+1-i}(-c_i-\varepsilon)$ itself multiplied by $(-1)^{3-k}$ (where $3-k$ is the positive inertia index of the critical point of $\check f$ with critical value $-c_i$). 
This proves the statement of our proposition about cycles vanishing in real critical points. The statement about the non-real ones is obvious.
(In fact, the exact value of the sign $\pm$ in it is also determined by the virtual Morse function of $f$, but its description is too long).
\hfill $\Box$

\begin{corollary}
\label{corud}
Under the conditions of Proposition \ref{proud}, the intersection indices of vanishing cycles satisfy the following relations:

1$)$ for any $i < j \leq r$, $\langle \Delta_i, \Delta_j \rangle = \langle \tilde \Delta_{r+1-i}, \tilde \Delta_{r+1-j} \rangle $;

2$)$ for any $i \in \{1, \dots, r\}$ and $j \in \{1, \dots, s\}$,  
$$\langle \Delta_i, \Delta_{r+j} \rangle = \pm \langle \tilde \Delta_{r+1-i}, \tilde \Delta_{r+s+1-j} \rangle \quad \mbox{ and } \quad 
\langle \Delta_i, \Delta_{r+s+j} \rangle = \pm \langle \tilde \Delta_{r+1-i} , \tilde \Delta_{9-j} \rangle.$$ 
\end{corollary}

\noindent
{\it Proof.} Both central symmetry and complex conjugation preserve the intersection forms in the homology of complex two-dimensional Milnor fibers of real polynomials. \hfill $\Box$

\subsection{D-graph invariant} 
In the case of polynomials whose critical points are all real, the set-valued invariant has the following more transparent interpretation.

Recall that the matrix of intersection indices $\langle \Delta_i, \Delta_j \rangle$ of basic vanishing cycles $\Delta_i , \Delta_j \in H_2(V_f)$ can be represented by the {\em Coxeter-Dynkin graph} (see e.g. \cite{AVGZ2}). Its vertices correspond to these vanishing cycles (in particular, they are numbered), and the $i$th and $j$th vertices are connected by $\langle \Delta_i, \Delta_j \rangle$ solid segments if this intersection index is positive, and by $-\langle \Delta_i, \Delta_j \rangle$ dashed segments if it is negative. 

\begin{definition} \rm
\label{dfdinv}The {\it D-graph} of a generic real Morse polynomial $f$ with only real critical points is \ $<$the isomorphism class of$>$ \ the oriented graph with vertices of two types (``even'' and ``odd''), obtained from the Coxeter-Dynkin graph of $f$ by

1) orienting each edge of this graph from the vertex corresponding to the critical point with the lower critical value to the one with the higher critical value;

2) indicating the parity of the Morse index of each critical point of $f$ at the corresponding vertex of the graph, and

3) forgetting the numbering of the vertices.
\end{definition} 

\noindent{\bf Notation.}
In future pictures of D-graphs (see Figs.~\ref{54}--\ref{1233} and \ref{324}--\ref{1413})
we will mark the ``odd'' (respectively, ``even'') vertices with white (respectively, black) little discs. \smallskip

Obviously, the D-graph of a Morse polynomial with only real critical points is determined by the associated virtual Morse function.

\begin{proposition}[cf. \cite{rigidplane}, Theorem 2]
\label{cont}
The D-graph is an invariant of isotopy classes of Morse polynomials of types $\Xi_1$ and $\Xi_2$ with eight real critical points. This invariant is equivalent to the set-valued invariant of Definition \ref{defsv}: any two such polynomials have equal values of one of these invariants if and only if they have equal values of the other. 
\end{proposition}

\noindent
{\it Proof} of this proposition consists of the following two lemmas.

\begin{lemma}
\label{p11}
Elementary virtual surgeries, which model the real surgeries and do not change the isotopy classes of Morse polynomials with only real critical points, do not change the corresponding D-graphs. In particular, the D-graph of a virtual Morse function is determined by its virtual component.
\end{lemma}

\noindent
{\it Proof}. These are only the surgeries $s2$, $s5$ and $s6$ from Definition \ref{elsur}. A surgery of type $s2$ does not change the D-graphs
because it preserves the set of homology classes of the vanishing cycles, only permuting two of them that are not connected by edges of the graph. Indeed, the cycles $\Delta_i(c_i - \varepsilon) \subset f^{-1}(c_i -\varepsilon)$ and $\Delta_{i+1}(c_{i+1} - \varepsilon) \subset f^{-1}(c_{i+1} -\varepsilon)$ are contained in small neighborhoods of distant critical points, therefore a continuous deformation of the function changing the order of its critical values does not affect the homology classes of vanishing cycles in $f^{-1}(0)$ obtained by transporting these cycles over continuously changing paths.

The surgery s5 or s6 can be realized as the composition of a) moving the basic non-critical value (in whose level set the vanishing cycles lie) along an arc starting at 0 and jumping over a neighboring critical value (while all paths to these non-critical values and the corresponding vanishing cycles in the fibers above them are continuously deformed, and the intersection indices do not change), and
\begin{figure}
\unitlength 0.7mm
\begin{picture}(48,20)
\put(0,5){\line(1,0){40}}
\bezier{150}(5,5)(10,15)(15,5)
\bezier{100}(20,5)(18,8)(15,5)
\bezier{250}(15,5)(25,20)(35,5)
\put(6,5){\circle*{1}}
\put(21,5){\circle*{1}}
\put(36,5){\circle*{1}}
\put(13.7,1.5){\footnotesize $0$}
\put(40,7){$\Longrightarrow$}
\end{picture} 
\begin{picture}(48,20)
\put(0,5){\line(1,0){40}}
\bezier{200}(5,5)(15,15)(20,10)
\put(20,5){\line(0,1){5}}
\bezier{200}(20,10)(25,15)(35,5)
\put(6,5){\circle*{1}}
\put(21,5){\circle*{1}}
\put(36,5){\circle*{1}}
\put(13.7,1.5){\footnotesize $0$}
\put(40,7){$\Longrightarrow$}
\end{picture} 
\begin{picture}(48,20)
\put(0,5){\line(1,0){40}}
\bezier{150}(5,5)(15,20)(25,5)
\bezier{100}(20,5)(22,8))(25,5)
\bezier{250}(25,5)(30,15)(35,5)
\put(6,5){\circle*{1}}
\put(21,5){\circle*{1}}
\put(36,5){\circle*{1}}
\put(13.7,1.5){\footnotesize $0$}
\put(40,7){$\Longrightarrow$}
\end{picture} 
\begin{picture}(50,20)
\put(0,5){\line(1,0){50}}
\bezier{150}(5,5)(15,20)(25,5)
\bezier{100}(20,5)(22,8))(25,5)
\bezier{250}(25,5)(30,15)(35,5)
\put(6,5){\circle*{1}}
\put(21,5){\circle*{1}}
\put(36,5){\circle*{1}}
\put(23.7,1.5){\footnotesize $0$}
\end{picture}
\caption{Decomposition of surgery $s6$}
\label{ss6}
\end{figure}
 b) adding a real constant to $f$, which moves this non-critical value back to its original position $0$, see Fig.~\ref{ss6}. These operations also preserve D-graphs. \hfill $\Box$ \smallskip

The orientations of the edges of the D-graph of any virtual Morse function with only real critical points define a partial order on the set of its vertices (since they obviously cannot form oriented cycles). \smallskip

\begin{lemma}[see Lemma 2 in \cite{rigidplane}]
\label{p111}
All virtual Morse functions defining the same D-graph can be obtained from each other by elementary surgeries $s2$, $s5$ and $s6$.

For any D-graph there is a natural one-to-one correspondence between the set of all virtual Morse functions defining it and the set of pairs consisting of 

a$)$ an isomorphism class of graphs with ordered vertices marked by the parity signs, which can be obtained by 
continuations of the partial order of vertices of our D-graph $($defined by the orientation of its edges$)$ to a complete order, and

b$)$ an integer number in the range from $0$ to $8$ $($describing the last element of the virtual Morse function, i.e., the number of negative critical values$)$. \hfill $\Box$ \end{lemma} 
\medskip

Proposition \ref{cont} follows immediately from Lemmas \ref{p11} and \ref{p111}. \hfill  $\Box$
\medskip

\begin{proposition}
\label{invofd}
The D-graphs of two Morse polynomials $f$ and $\check f$ related with each other by the ``up-down'' involution $($\ref{invof}$)$ are obtained from each other by reversing the orientations of all edges and changing the colors of all vertices.
\end{proposition}

\noindent
{\it Proof.} Let us associate the vertex of the D-graph of $f$ corresponding to the $i$th critical value with the vertex of the D-graph of $\check f$ corresponding to the $(9-i)$th critical value. This association obviously changes all the parities of the Morse indices and reverses the orientation of all the edges. By Corollary \ref{corud} it also preserves the numbers of edges connecting the associated pairs of vertices. \hfill $\Box$

\subsection{Independent perturbations of critical points} 
\label{indper}

A neighborhood of a polynomial of degree three with a non-discriminant principal homogeneous part in the space of all such polynomials can be considered as a piece of the parameter space of a versal deformation of the polynomial given by this principal part.

 Therefore, for any polynomial of this type with several critical points, any collection of small perturbations of its restrictions on small neighborhoods of these points can be realized (up to a set of coordinate changes in these neighborhoods) by a common perturbation of our polynomial in the corresponding class $\Xi_1$ or $\Xi_2$.

\section{Fundamental lemma}

\begin{lemma}
A third degree polynomial $f: {\mathbb R}^3 \to {\mathbb R}$ with 
only isolated critical points in ${\mathbb C}^3$ cannot have four Morse real critical points lying in the same affine plane and having the same parities of the Morse indices.
\label{flem}
\end{lemma}

\noindent
{\it Proof.} 
By Proposition \ref{no2max}, there are only four possibilities for the Morse indices of such four critical points:
\begin{enumerate}
\item one point of index 3 and three points of index 1,
\item \label{tww} four points of index 1,
\item \label{trr} one point of index 0 and three points of index 2,
\item four points of index 2.
\end{enumerate}
We can and will consider only the cases (\ref{tww}) and (\ref{trr}), because the other two 
can be reduced to them by multiplying the function by $-1$.

Suppose that $f$ is a polynomial whose four Morse critical points with one of these two index collections lie in the plane $\{z=0\}$. Let
\begin{equation}
f (x, y, z) \equiv \varphi(x,y) + z f_2(x,y) + z^2 f_1(x,y) + f_0 z^3 
\label{decom}
\end{equation}
be its decomposition by the powers of $z$, where $\varphi$, $f_2$ and $f_1$ are polynomials in the coordinates $x$ and $ y$ of degree at most $3, 2,$ and $1$, respectively. Obviously, all critical points of $f$ lying in the plane ${\mathbb R}^2 = \{z=0\}$ are also critical points of $\varphi$.

\begin{proposition}
The homogeneous principal part of the polynomial $\varphi$ is not dis\-cri\-minant $($i.e., it vanishes on three different lines in ${\mathbb C}^2)$. 
\end{proposition}

\noindent
{\it Proof.} If this principal part is discriminant but $\varphi$ has only isolated critical points in ${\mathbb C}^2$, then by the Bezout's theorem the number of these points is less than 4. If $\varphi$ has non-isolated critical points but $\varphi \not \equiv 0$ then all its real critical points form a line in ${\mathbb R}^2$. The critical points of the entire function $f$ lying in ${\mathbb R}^2$ are the intersection points of this line with the zero set of the quadratic polynomial $f_2$. This intersection set is either infinite (in which case the entire function $f$ has infinitely many critical points) or consists of at most two points. Finally, if $\varphi \equiv 0$ then the polynomial $f$ is reducible and has non-isolated critical points in ${\mathbb C}^3$. \hfill $\Box$ \smallskip

\begin{proposition}
\label{pro71}
If a polynomial $\varphi:{\mathbb R}^2 \to {\mathbb R}$ of degree three with non-discriminant principal homogeneous part has four real Morse critical points, then one of two things is true:
\begin{enumerate}
\item 
the principal homogeneous part of $\varphi$ vanishes on only one real line, and the critical points of $\varphi$ in ${\mathbb R}^2$ are a minimum, a maximum and two saddlepoints; the convex hull of these four points is a quadrilateral with one diagonal connecting two local extrema and the other diagonal connecting two saddlepoints;
\item
the principal homogeneous part of $\varphi$ vanishes on three real lines, and the critical points of $\varphi$ are three saddlepoints and one point of local extremum, which lies strictly inside the triangle with vertices at the saddlepoints.
\end{enumerate}
\label{pro19}
\end{proposition}

\noindent
{\it Proof.} By an affine transformation of ${\mathbb R}^2$ and adding a constant function we can reduce $\varphi$ to the form 
\begin{equation}
x^2 y \pm y^3 + a x + b y + c y^2, 
\label{d4}
\end{equation} 
which is actually the normal form of a versal deformation of $D_4^+$ or $D_4^-$ singularity reduced modulo the addition of the constants, cf. \cite{AVGZ1}.
\begin{figure}
\unitlength 0.8mm
\linethickness{0.4pt}
\begin{center}
\begin{picture}(47.00,54.00)
\bezier{264}(0.00,17.00)(15.00,-13.00)(37.00,12.00)
\bezier{44}(12.00,16.00)(17.00,12.00)(17.00,7.00)
\bezier{52}(17.00,7.00)(19.00,13.00)(25.00,15.00)
\bezier{164}(37.00,12.00)(32.00,32.00)(36.00,52.00)
\bezier{176}(0.00,17.00)(15.00,34.00)(23.00,54.00)
\bezier{152}(12.00,16.00)(13.00,41.00)(11.00,54.00)
\bezier{168}(25.00,15.00)(38.00,37.00)(47.00,51.00)
\bezier{52}(26.00,42.00)(26.00,50.00)(23.00,54.00)
\bezier{60}(26.00,42.00)(29.00,49.00)(36.00,52.00)
\put(20.10,26.00){\circle{1.00}}
\bezier{52}(19.70,26.10)(16.00,29.50)(12.20,33.20)
\bezier{64}(20.40,26.20)(27.10,29.15)(34.20,30.35)
\bezier{68}(20.30,26.40)(23.00,37.00)(26.00,42.00)
\bezier{76}(20.00,25.50)(18.00,17.00)(17.00,7.00)
\bezier{52}(17.00,43.00)(14.00,47.00)(11.00,54.00)
\bezier{64}(35.00,40.00)(40.00,42.50)(47.00,51.00)
\bezier{32}(33.00,39.00)(29.00,37.00)(26.00,38.00)
\bezier{16}(19.00,41.00)(20.00,40.00)(22.00,39.00)
\end{picture}
\qquad \qquad \qquad \quad
\begin{picture}(40.00,58.00)
\put(22.00,27.00){\circle{1.00}}
\bezier{84}(22.30,27.40)(27.00,35.00)(23.00,46.00)
\bezier{124}(22.30,27.40)(29.00,38.00)(31.00,55.00)
\bezier{96}(22.30,27.40)(29.00,38.00)(38.00,45.00)
\bezier{56}(23.00,46.00)(28.00,47.00)(31.00,55.00)
\bezier{56}(38.00,45.00)(33.00,47.00)(31.00,55.00)
\bezier{84}(21.70,26.60)(19.00,20.00)(7.00,15.00)
\bezier{84}(7.00,15.00)(11.00,8.00)(24.00,7.00)
\bezier{80}(24.00,7.00)(13.00,7.00)(5.00,4.00)
\bezier{48}(5.00,4.00)(7.00,7.00)(7.00,15.00)
\bezier{36}(5.00,4.00)(9.00,6.00)(12.50,9.60)
\bezier{96}(24.00,7.00)(17.00,14.00)(21.70,26.60)
\bezier{12}(25.00,46.00)(27.00,46.00)(28.00,45.90)
\bezier{16}(31.00,45.85)(32.00,45.84)(35.50,45.50)
\bezier{30}(13.7,11.3)(17,15)(19,20)
\end{picture}
\caption{Caustics for degree 3 polynomials in two variables: {\em purse} ($D_4^+$) and {\em pyramid} ($D_4^-$)}
\label{pyramid}
\end{center}
\end{figure}
 The variety of non-Morse functions in the space of such polynomials with principal part $x^2y+y^3$ (i.e., of class $D_4^+$) is known (see, for example, \S~I.3 in \cite{AVGZ1} or \S~ II.2 in \cite{AGLV1}). Its intersection with a neighborhood of the origin in the parameter space of the family (\ref{d4}) is shown in Fig.~\ref{pyramid} (left). The polynomials with four real critical points occupy there one of the four connected components of the complement of this variety (the upper component in this picture), and have the sets of Morse indices described in point 1 of Proposition \ref{pro71}.

In the space of polynomials (\ref{d4}) with principal part $x^2 y - y^3$ (of class $D_4^-$) the polynomials with four critical points lie inside two pyramids of Fig.~\ref{pyramid} (right) and have three saddlepoints and a local minimum or maximum depending on the pyramid. 

It is easy to see that a polynomial of the form (\ref{d4}) cannot have three critical points on the same affine line, so by the connectivity of our components the geometric disposition of four critical points in the plane either satisfies or does not satisfy the conditions of Proposition \ref{pro71} for all parameters $a, b, c$ within these components. 

In both cases $D_4^+$ and $D_4^-$, it is easy to construct the polynomials with the sets of critical points that satisfy the corresponding condition, hence all these polynomials satisfy them. \hfill $\Box$ 

\begin{corollary}
\label{cora}
1. If all four critical points of $f$ lying in the plane $\{z=0\}$ have Morse indices equal to 1, then the possibility $($2$)$ of Proposition \ref{pro71} holds for the corresponding critical points of $\varphi$; moreover, the local extremum of $\varphi$ is a local minimum. 

2. If one of these four critical points of $f$ has Morse index equal to 0 and the other three have Morse indices equal to 2, then either the possibility $($1$)$ of Proposition \ref{pro71} holds or the special case of the possibility $($2$)$ holds, in which $\varphi$ has three saddlepoints and one local minimum.
\end{corollary}

\noindent
{\it Proof.} 1. If the restriction $\varphi$ of the function $f$ to the plane ${\mathbb R}^2$ has a local maximum, then the Morse index of $f$ at that point is at least 2. 

2. The local minimum of $f$ is also a local minimum of $\varphi$.
\hfill $\Box$
\smallskip


The condition that our four critical points of $\varphi$ are also critical points of $f$ implies that the polynomial $f_2$ of (\ref{decom}) vanishes at these points. The polynomials ${\mathbb R}^2 \to {\mathbb R}$ of degree 2 that vanish at these four points form a two-dimensional vector space. The polynomials $\varphi'_x$ and $\varphi'_y$ are linearly independent in this space because $\varphi$ has only isolated singularities. Therefore, 
\begin{equation}
f_2 \equiv \alpha \varphi'_x + \beta \varphi'_y
\end{equation}
for some constants $\alpha$ and $\beta$. The quadratic part of the Taylor expansion of the polynomial $f$ at a critical point $K$ of $\varphi$ is then given by the matrix 
$$
\frac{1}{2}\left|\mbox{
\begin{tabular}{lll}
$\varphi''_{x x}(K)$ & $\varphi''_{x y}(K)$ & $\alpha \varphi''_{x x}(K) + \beta \varphi''_{x y}(K)$ \\
$\varphi''_{y x}(K)$ & $\varphi''_{y y}(K)$ & $\alpha \varphi''_{x y}(K) + \beta \varphi''_{y y}(K)$ \\
$\alpha \varphi''_{x x}(K) + \beta \varphi''_{x y}(K)$ & $\alpha \varphi''_{x y}(K) + \beta \varphi''_{y y}(K)$ & $ 2 f_1(K)$
\end{tabular}
}\right|.
$$
The determinant of this quadratic form is a linear function of the value $f_1(K)$ and vanishes if and only if $2 f_1(K)$ is equal to the value of the linear function 
\begin{equation}\label{eee} 
\alpha^2 \varphi''_{x x} + 2 \alpha \beta \varphi''_{x y} + \beta^2 \varphi''_{y y}
\end{equation}
at the point $K$.
If $2 f_1(K)$ is greater (respectively, less) than this value, then the signature of the Morse critical point $K$ of $f$
(i.e., the pair $(i_+, i_-)$ of positive and negative inertia indices of its quadratic part) is equal to the sum of 
a) the signature of the same point considered as a critical point of $\varphi$, and b) the string $(1,0)$ (respectively, $(0,1)$).

In the situation considered in point 1 of Corollary \ref{cora}, all our four critical points of $f$ are of signature $(2, 1)$, which implies that the linear function $2 f_1$ is greater than the linear function 
(\ref{eee}) 
at all three corners of a triangle, but is less than it at an interior point of the same triangle. In the situation of point 2 of Corollary \ref{cora}, the signatures of four points of $f$ are equal to (3,0), (1,2), (1,2), (1,2), and the signatures of the corresponding points of $\varphi$ are either $(2,0), (1,1),$ $(1,1), (1,1)$ or $(2,0), (1,1),$ $(1,1), (0,2)$. In the first case the linear function $2 f_1$ is less than (\ref{eee}) at all the corners of a triangle, but greater than it at some of its interior points. In the second case, $2 f_1$ is greater than (\ref{eee}) at the endpoints of one diagonal of a convex plane quadrilateral but less than it at the endpoints of the other diagonal. All these systems of inequalities are impossible, so we have a contradiction.
\hfill $\Box$ 

\begin{remark} \rm
All combinations of four Morse indices not forbidden by Lemma \ref{flem} and Proposition \ref{no2max} can be realized by critical points of polynomials of type $\Xi_1$ or $\Xi_2$ lying in the same affine plane.
\end{remark}

\section{Isotopy classes of Morse polynomials of type $\Xi_1$}
\label{p81}

\subsection{Polynomials with eight critical points}
\label{p8eight}

\begin{proposition}
\label{pro14}
The set of all 3789 virtual Morse functions of type $\Xi_1$ with eight real critical points consists of exactly seven virtual components: 
\begin{itemize}
\item one of cardinality 54 with the D-graph shown in Fig.~\ref{54};

\item one of cardinality 72 with the D-graph shown in Fig.~\ref{72};

\item two of cardinality 432: the D-graph of one of them is shown in Fig.~\ref{432}, and that of the other is obtained from it by the involution from Proposition \ref{invofd};

\item two of cardinality 783: the D-graph of one of them is shown in Fig.~\ref{783}, and that of the other is obtained from it by the involution from Proposition \ref{invofd};

\item one of cardinality 1233 with D-graph shown in Fig.~\ref{1233}.
\end{itemize} 
\end{proposition}

\begin{figure}
\unitlength 0.7mm
\begin{center}
\begin{picture}(123,53)
\put(0,5){\circle{1.9}} 
\put(40,5){\circle{1.9}}
\put(80,5){\circle{1.9}}
\put(120,5){\circle{1.9}}
\put(0,45){\circle*{1.9}}
\put(40,45){\circle*{1.9}}
\put(80,45){\circle*{1.9}}
\put(120,45){\circle*{1.9}}
\put(0,36){\vector(0,1){7}}
\put(0,34){\line(0,-1){8}}
\put(0,24){\line(0,-1){8}}
\put(0,14){\line(0,-1){8}}
\put(40,36){\vector(0,1){7}}
\put(40,34){\line(0,-1){8}}
\put(40,24){\line(0,-1){8}}
\put(40,14){\line(0,-1){8}}
\put(120,36){\vector(0,1){7}}
\put(120,34){\line(0,-1){8}}
\put(120,24){\line(0,-1){8}}
\put(120,14){\line(0,-1){8}}
\put(79.6,14){\vector(0,-1){7}}
\put(80.4,14){\vector(0,-1){7}}
\put(79.6,16){\line(0,1){8}}
\put(80.4,16){\line(0,1){8}}
\put(79.6,26){\line(0,1){8}}
\put(80.4,26){\line(0,1){8}}
\put(79.6,36){\line(0,1){8}}
\put(80.4,36){\line(0,1){8}}
\put(0.8,5.4){\vector(2,1){78}}
\put(40.7,5.7){\vector(1,1){38.5}}
\put(119.3,5.7){\vector(-1,1){38.5}}
\put(100.5,25.5){\vector(1,1){19}}
\put(99.5,24.5){\line(-1,-1){19}}
\put(39.5,25.25){\vector(-2,1){39.2}}
\put(40.5,24.7){\line(2,-1){12}}
\put(53.6,18.05){\line(2,-1){25.3}}
\put(52.9,32.1){\vector(-1,1){12.3}}
\bezier{100}(53.8,31.2)(55.8,29.2)(59.6,25.4)
\put(60.5,24.5){\line(1,-1){18.8}}
\put(81,45){\vector(1,0){38}}
\put(79,45){\vector(-1,0){38}}
\put(119,5){\vector(-1,0){38}}
\put(41,5){\vector(1,0){38}}
\bezier{600}(0.7,4.3)(60,-5)(78,4)
\put(78,4){\vector(2,1){1}}
\bezier{600}(79.3,45.7)(20,55)(1,46)
\put(1,46){\vector(-2,-1){0.5}}
\end{picture}
\end{center}
\caption{D-graph related to bisingularity $D_4^- + D_4^-$; Card = 54}
\label{54}
\end{figure}

\begin{figure}
\unitlength 0.7mm
\begin{center}
\begin{picture}(90,85)
\put(0,5){\circle*{1.9}} 
\put(0,65){\circle{1.9}} 
\put(60,5){\circle{1.9}} 
\put(30,20){\circle{1.9}}
\put(0,6){\vector(0,1){58}}
\put(1,5){\vector(1,0){58}}
\put(0.8,5.6){\vector(2,1){28.3}}
\put(30,80){\circle*{1.9}}
\put(60,65){\circle*{1.9}}
\put(90,20){\circle*{1.9}}
\put(90,80){\circle{1.9}}
\put(1,65){\vector(1,0){58.3}}
\put(60,6){\vector(0,1){58.3}}
\put(30,71.4){\vector(0,1){8}}
\put(30,69.6){\line(0,-1){4.1}}
\put(30,64.4){\line(0,-1){28.4}}
\bezier{70}(30,34.4)(30,32.5)(30,30.7)
\put(30,29.4){\line(0,-1){8.7}}
\put(0.8,65.6){\vector(2,1){28.3}}
\put(60.5,20){\vector(1,0){28.5}}
\put(59.5,20){\line(-1,0){28.7}}
\put(60.8,5.6){\vector(2,1){28.4}}
\put(90,21){\vector(0,1){58}}
\put(31,80){\vector(1,0){57.8}}
\put(60.8,65.6){\vector(2,1){28.4}}
\put(24.6,66.5){\line(2,5){5}}
\put(29,77.5){\vector(1,2){0.5}}
\put(22,60){\line(-2,-5){5.6}}
\put(14,40){\line(-2,-5){5.6}}
\put(6,20){\line(-2,-5){5.4}}
\put(89.3,79.3){\line(-2,-5){5.6}}
\put(89,77.5){\vector(1,2){0.5}}
\put(82,60){\line(-2,-5){5.6}}
\put(74,40){\line(-2,-5){5.6}}
\put(65.4,18){\line(-2,-5){4.8}}
\put(51,56){\vector(1,1){8.2}}
\put(47,52){\line(-1,-1){8}}
\put(35,40){\line(-1,-1){8}}
\put(23,28){\line(-1,-1){8}}
\put(11,16){\line(-1,-1){10}}
\put(81,71){\vector(1,1){8}}
\put(77,67){\line(-1,-1){8}}
\bezier{120}(65,55)(63,53)(60.5,50.5)
\bezier{130}(59.5,49.5)(57.7,47.7)(55.7,45.7)
\put(52,42){\line(-1,-1){8}}
\put(41,31){\line(-1,-1){10.5}}
\put(90,20){\line(-6,-1){12}}
\put(87,19.5){\vector(4,1){0.5}}
\put(73.8,17.3){\line(-6,-1){8.7}}
\bezier{50}(63.9,15.65)(62.7,15.45)(60.9,15.15)
\put(57.6,14.6){\line(-6,-1){12}}
\put(41.4,11.9){\line(-6,-1){12}}
\put(25.2,9.2){\line(-6,-1){12}}
\put(9,6.5){\line(-6,-1){8}}
\put(89.4,79.9){\line(-6,-1){11.4}}
\put(87,79.5){\vector(4,1){0.5}}
\put(73.8,77.3){\line(-6,-1){12}}
\put(57.6,74.6){\line(-6,-1){12}}
\put(41.4,71.9){\line(-6,-1){12.2}}
\put(25.2,69.2){\line(-6,-1){12}}
\put(9,66.5){\line(-6,-1){8}}
\put(89.4,79.5){\line(-6,-5){28.8}}
\put(59.6,54.7){\line(-6,-5){59.5}}
\end{picture}
\end{center}
\caption{$D_4^+ + D_4^+$, $D_5+ A_3$, and $D_4^- + D_4^-$; Card = 72}
\label{72}
\end{figure}

\begin{figure}
\unitlength 0.7mm
\begin{center}
\begin{picture}(123,88)
\put(0,5){\circle{1.9}} 
\put(40,5){\circle{1.9}}
\put(80,5){\circle{1.9}}
\put(0,45){\circle*{1.9}}
\put(40,45){\circle*{1.9}}
\put(80,45){\circle*{1.9}}
\put(120,45){\circle*{1.9}}
\put(60,85){\circle{1.9}} 
\put(0,5.8){\vector(0,1){38.5}}
\put(39.2,5.6){\vector(-1,1){38.5}}
\put(79.2,5.4){\vector(-2,1){78.5}}
\put(0.7,45.7){\vector(3,2){58.3}}
\put(40.4,45.7){\vector(1,2){18.8}}
\put(79.6,45.7){\vector(-1,2){18.8}}
\put(119.2,45.6){\vector(-3,2){58.3}}
\put(27,32){\vector(1,1){12.3}}
\bezier{100}(26.1,31.1)(23.1,28.1)(20.4,25.4)
\put(19.3,24.3){\line(-1,-1){18.8}}
\put(53.9,18.9){\vector(1,1){25.2}}
\put(52.8,17.8){\line(-1,-1){12.3}}
\put(80.7,5.7){\vector(1,1){38.4}}
\put(49.2,70.6){\vector(3,4){10}}
\put(46.5,67){\line(-3,-4){10.5}}
\put(33,49){\line(-3,-4){10}}
\bezier{60}(21,33)(19.5,31)(18,29)
\put(14.1,23.8){\line(-3,-4){13.6}}
\put(56.65,71.6){\vector(1,4){3.25}}
\put(55.75,68){\line(-1,-4){3.25}}
\put(51.5,51){\line(-1,-4){3.25}}
\put(47.25,34){\line(-1,-4){2.7}}
\put(43.5,19){\line(-1,-4){3.3}}
\put(63.35,71.6){\vector(-1,4){3.25}}
\put(64.25,68){\line(1,-4){3.25}}
\put(68.5,51){\line(1,-4){3.25}}
\put(72.75,34){\line(1,-4){3.25}}
\put(77,17){\line(1,-4){2.75}}
\end{picture}
\end{center}
\caption{$D_4^- + D_4^-$, $E_6 + A_2$, and $D_5+A_3$; Card = 432}
\label{432}
\end{figure}

\begin{figure}
\unitlength 0.7mm
\begin{center}
\begin{picture}(123,48)
\put(0,5){\circle{1.9}} 
\put(40,5){\circle{1.9}}
\put(120,5){\circle{1.9}}
\put(120,45){\circle*{1.9}}
\put(0,45){\circle*{1.9}}
\put(40,45){\circle*{1.9}}
\put(80,45){\circle*{1.9}}
\put(80,5){\circle{1.9}} 
\put(0,6){\vector(0,1){38.5}}
\put(0.8,5.6){\vector(2,1){78.5}}
\put(40,26){\vector(0,1){18}}
\put(40,24.6){\line(0,-1){18.5}}
\put(40.7,5.7){\vector(1,1){38.5}}
\put(119.3,5.7){\vector(-1,1){38.5}}
\put(120,8){\line(0,1){7}}
\put(120,18){\line(0,1){7}}
\put(120,28){\line(0,1){7}}
\put(120,38){\vector(0,1){6}}
\put(41,5){\vector(1,0){38}}
\put(119,5){\vector(-1,0){38}}
\bezier{600}(0.7,4.3)(50,-5)(79,4)
\put(80.7,5.7){\line(1,1){18.5}}
\put(100.7,25.7){\vector(1,1){18}}
\put(119,45){\vector(-1,0){38}}
\put(79,4){\vector(2,1){0.4}}
\put(80.4,14){\vector(0,-1){8}}
\put(79.6,14){\vector(0,-1){8}}
\put(80.4,16){\line(0,1){8}}
\put(79.6,16){\line(0,1){8}}
\put(80.4,26){\line(0,1){8}}
\put(79.6,26){\line(0,1){8}}
\put(80.4,36){\line(0,1){8}}
\put(79.6,36){\line(0,1){8}}
\end{picture}
\end{center}
\caption{$E_6+A_2$, Card = 783}
\label{783}
\end{figure}

\begin{figure}
\unitlength 0.7mm
\begin{center}
\begin{picture}(123,53)
\put(0,5){\circle{1.9}} 
\put(40,5){\circle{1.9}}
\put(80,5){\circle{1.9}}
\put(120,5){\circle{1.9}}
\put(0,45){\circle*{1.9}}
\put(40,45){\circle*{1.9}}
\put(80,45){\circle*{1.9}}
\put(120,45){\circle*{1.9}}
\put(0,6){\vector(0,1){38.3}}
\put(40,6){\vector(0,1){38.3}}
\put(80,6){\vector(0,1){38.3}}
\put(120,6){\vector(0,1){38.3}}
\put(0.7,5.7){\vector(1,1){38.6}}
\put(40.8,45){\vector(1,0){38}}
\put(41,5){\vector(1,0){38.3}}
\put(80.7,5.7){\vector(1,1){38.5}}
\bezier{600}(0.7,4.5)(60,-5)(79,4)
\put(79,4){\vector(2,1){0.5}}
\bezier{600}(40.7,45.5)(100,55)(119,46)
\put(119,46){\vector(2,-1){0.5}}
\put(71,36){\vector(1,1){8}}
\put(69,34){\line(-1,-1){8}}
\put(59,24){\line(-1,-1){8}}
\put(49,14){\line(-1,-1){8}}
\put(71,14.5){\vector(1,-1){8}}
\put(70,13.5){\vector(1,-1){8}}
\put(69,16.5){\line(-1,1){8}}
\put(68,15.5){\line(-1,1){8}}
\put(59,26.5){\line(-1,1){8}}
\put(58,25.5){\line(-1,1){8}}
\put(49,36.5){\line(-1,1){7.5}}
\put(48,35.5){\line(-1,1){7.5}}
\end{picture}
\end{center}
\caption{$A_4+A_4$, Card = 1233}
\label{1233}
\end{figure}

\noindent
{\it Proof.}
A slightly modified version 
{\footnotesize
\begin{verbatim}
https://drive.google.com/file/d/1TAx9X9Rhgt6A-ALizF5KG_gBlOg6MxpO/view?usp=sharing
\end{verbatim}
\label{pro22}}
\noindent
of the program used in the proof of Proposition \ref{procount1} (with virtual surgeries of types $s1$ and $s3$ disabled) enumerates the elements of the virtual component of our initial virtual Morse function (with the D-graph shown in Fig.~\ref{72}). The number of these elements is smaller than the number of all virtual functions of type $\Xi_1$ with eight real critical points, so we choose another such virtual Morse function found by our first program, check if its D-graph is different from the previous one (if not, we simply repeat the search in another sector of the formal graph), and compute its $\mbox{Card}$ value again. So we get another virtual component (or even two of them if the D-graph of the new component is not invariant with respect to the up-down involution from Proposition \ref{invofd}) and continue the search. The procedure stops when the sum of the Card invariants of all found virtual components is equal to the total number 3789 of virtual Morse functions of type $\Xi_1$ with eight real critical points. \hfill $\Box$ \smallskip  

\begin{remark} \rm
This modified program can generally start with virtual Morse functions, whose critical points are not all real. For this reason, the second and third lines from the bottom of the initial data file like (\ref{virtuM}) (right) must also be inserted as values $\mbox{PC1}(i)$ and $\mbox{PC2}(i)$, $i=1, \dots, 8$, in the last part of this program. Also, the number MIC of pairs of non-real critical points should be specified in line 37.
\end{remark}

\begin{proposition}
\label{protrue}
All white $($respectively, black$)$ vertices in figures \ref{54}, \ref{783}, and \ref{1233} correspond to critical points with Morse index 1 $($respectively, 2$)$ of all Morse polynomials with these D-graphs. The white vertex at the top of Fig.~\ref{432} corresponds to a point of local maximum, and all other white $($respectively, black$)$ vertices of this graph correspond to points with Morse index 1 $($respectively, 2$)$.
The top right white vertex in Fig.~\ref{72} corresponds to a local maximum, the bottom left black vertex corresponds to a local minimum, and all other white $($respectively, black$)$ vertices correspond to critical points with Morse index 1 $($respectively, 2$)$.
\end{proposition}

\noindent
{\it Proof.} This follows easily from Proposition \ref{lemfun}, cf. Example on p. \pageref{exaa}.
\hfill $\Box$

\begin{theorem}
\label{achi}
1. All virtual components of type $\Xi_1$ with eight real critical points, except for the virtual component with Card=1233 and D-graph shown in Fig.~\ref{1233}, are achiral and thus are associated with single isotopy classes of Morse polynomials.

2. The virtual component with Card=1233 is chiral and thus associated with two isotopy classes which are mapped into each other by the composition of polynomials with the reflection in any plane in ${\mathbb R}^3$.

In particular, there are exactly eight isotopy classes of degree three Morse polynomials ${\mathbb R}^3 \to {\mathbb R}$ of type $\Xi_1$ with eight real critical points. 
\end{theorem}

The proof of this theorem is given in \S~\ref{prachi}.

\subsubsection{Realization of Figs.~\ref{54}--\ref{432}}

The homogeneous polynomial 
\begin{equation} 
\label{mex}
x^3 - 3 x y^2 + z^3 
\end{equation}
is of type $\Xi_1$. Indeed, by adding $\varepsilon z$, $\varepsilon >0$, we remove all its real critical points, which is impossible for polynomials of type $\Xi_2$, see \S~\ref{morcom}. Also, it is easy to check that it has an isolated singularity at the origin.

Another its perturbation,
\begin{equation}
\label{zero}
x^3 - 3 x y^2 + z^3 - 3 \varepsilon z , \qquad \varepsilon >0,
\end{equation}
has two critical points of the class $D_4^-$ (see Table \ref{t1} on page \pageref{t1}). In a neighborhood of one of them this polynomial is equivalent to 
\begin{equation}
\label{one}
 x^3 - 3 x y^2 + \tilde z^2 - 2 \varepsilon^{3/2} ,
\end{equation}
and in a neighborhood of the other one to 
\begin{equation}
\label{two}
 x^3 - 3 x y^2 - \tilde z^2 + 2 \varepsilon^{3/2} .
\end{equation}

Each of these two critical points can be further split into four Morse critical points by a small perturbation of the polynomial (\ref{zero}). 
Namely, the polynomial 
\begin{equation}x^3 - 3 x y^2 + z^3 - 3 \varepsilon z - \delta z (x^2 + y^2)
\label{twop}
\end{equation}
with sufficiently small $\delta>0$ has four points of signature $(1,2)$ and four points of signature $(2,1)$.
The polynomial 
\begin{equation}
x^3 - 3 x y^2 + z^3 - 3 \varepsilon z + \delta z (x^2 + y^2)
\label{onep}
\end{equation}
has one Morse point of signature $(0, 3)$, three points of signature $(1,2)$, three points of signature $(2,1)$ and one point of signature $(3,0)$.
The polynomial 
\begin{equation}x^3 - 3 x y^2 + z^3 - 3 \varepsilon z - \delta (x^2 + y^2)
\label{threep}
\end{equation}
has three critical points of signature $(2,1)$, four points of signature $(1,2)$ and one point of signature $(0,3)$. Its up-down version, 
\begin{equation}x^3 - 3 x y^2 + z^3 - 3 \varepsilon z + \delta (x^2 + y^2),
\label{fourp}
\end{equation}
has one point of signature $(3,0)$, four points of signature $(2,1)$ and three points of signature $(1,2)$. So, the passports of the Morse polynomials (\ref{twop}), (\ref{onep}), (\ref{threep}), and (\ref{fourp}) are $(0, 4, 4, 0)$, $(1, 3, 3, 1)$, $(0, 3, 4, 1)$, and $(1, 4, 3, 0),$ respectively. 

\begin{proposition}
The D-graphs of these four polynomials are respectively the D-graphs shown in Figs.~\ref{54}, \ref{72}, \ref{432}, and the graph obtained from Fig.~\ref{432} by the up-down involution described in Proposition \ref{invofd}. 
\end{proposition}

\noindent
{\it Proof.}
In all cases, the matrix of the intersection indices of four vanishing cycles arising from either of the critical points (\ref{one}) or (\ref{two}) is represented by the standard Coxeter--Dynkin graph 
\unitlength 0.4 mm
\begin{picture}(29,12)
\put(3,4){\circle*{2}}
\put(15,4){\circle*{2}}
\put(25,-2){\circle*{2}}
\put(25,10){\circle*{2}}
\put(3,4){\line(1,0){12}}
\put(15,4){\line(5,3){10}}
\put(15,4){\line(5,-3){10}}
\end{picture} 
of type $D_4$.
Therefore, for each of the four polynomials (\ref{twop})--(\ref{fourp}), we can divide the set of vertices of the corresponding D-graph into two subsets of cardinality four in such a way that, by removing all edges connecting vertices from different subsets, we split the D-graph into two such graphs 
(with some orientation of the edges and coloring of the vertices depending on the choice of perturbations). The D-graphs of Figs.~\ref{783} and \ref{1233} and their images under the up-down involution do not have such splittings. Thus, according to Proposition \ref{pro14}, the D-graphs of the polynomials (\ref{twop})--(\ref{fourp}) are among the remaining four D-graphs described in this proposition.

For each of the D-graphs of Figs.~\ref{54}, \ref{72}, \ref{432} and the up-down reflection of Fig.~\ref{432},
 the Morse indices of the critical points of all polynomials having these D-graphs
are determined by Proposition \ref{protrue}. In particular, the passports of these polynomials are equal to  $(0, 4, 4, 0)$, $(1, 3, 3, 1)$, $(0, 3, 4, 1)$, and $(1, 4, 3, 0)$ respectively. This fixes the only possible correspondence between the polynomials (\ref{twop})--(\ref{fourp}) and these four D-graphs.
 \hfill $\Box$

\subsubsection{Realization of Fig.~\ref{783}}
\label{rea783}
The polynomial 
\begin{equation}
\label{e66}
\frac{1}{3} x^3 + \frac{1}{3}z^3 - y^2 z + z^2 
\end{equation} 
(which is obviously also of type $\Xi_1$) 
has a critical point of class $E_6$ at the origin and is there equivalent to the function germ $ \frac{1}{3}x^3+\tilde z^2 - \frac{1}{4} y^4$. It also has the critical point $\left(0, 0, -2\right)$ of class $A_2$, where it is locally equivalent to the function $\frac{1}{3}\check x^3 + 2 \check y^2 - \check z^2 + \frac{4}{3}$. These two critical points can be simultaneously perturbed as follows: the point of class $A_2$ splits into two points with Morse indices 1 and 2, and the point of class $E_6$ splits into three Morse points with Morse index 1 and three points with Morse index 2 (see Fig.~\ref{e6proof} and p. 16 in \cite{AC} for the zero set of the corresponding perturbation of an $E_6$ singularity in two variables). Namely, such a perturbation is realized by the polynomial 
\unitlength 0.6mm 
\begin{figure} 
\begin{center}
{
\begin{picture}(50,28)
\bezier{300}(15,1)(40,31)(47,26)
\bezier{70}(47,26)(50,24)(47,22)
\bezier{250}(47,22)(25,14)(3,22)
\bezier{70}(3,22)(0,24)(3,26)
\bezier{300}(3,26)(10,31)(35,1)
\end{picture}
}\end{center}
\caption{Perturbation of $E_6$ singularity}
\label{e6proof}
\end{figure}
\begin{equation}
\label{reali77}
\frac{1}{3}x^3 + \frac{1}{3}z^3 - y^2 z + z^2 - \varepsilon x y^2 - \frac{\varepsilon^4}{12}x - \frac{\varepsilon^3}{2}y^2 
\end{equation}
with sufficiently small $\varepsilon>0$, 
cf. \cite{GZ}.

The corresponding D-graph should contain two subgraphs with complementary sets of 2 and 6 vertices, which are isomorphic to classical Coxeter--Dynkin graphs of types $A_2$ and $E_6$, and all edges of this D-graph connecting these subgraphs are directed from the $E_6 $ graph to the $A_2$ graph. Of all the D-graphs mentioned in Proposition \ref{pro14} (including the up-down versions of the graphs explicitly drawn there), only the one from Fig.~\ref{783} satisfies this property. So, the polynomial (\ref{reali77}) realizes exactly this D-graph, and minus this polynomial realizes its up-down version.

\begin{remark} \rm
\label{e67r}
The polynomial 
\begin{equation}
\label{e67}
\frac{1}{3} x^3 - \frac{1}{3}z^3 - y^2 z + z^2 
\end{equation} 
also has two critical points of classes $E_6$ and $A_2$. However, this polynomial is locally equivalent to $\frac{1}{3} x^3 - 2 y^2 - \check z^2 + \frac{4}{3}$ at its point of class $A_2$, and so any of its small perturbations with eight real critical points will have a local maximum.
\end{remark}

\subsubsection{Realization of Fig.~\ref{1233}.} 
The  polynomial \begin{equation}
\frac{20}{3} y^3 + x^2z + y z^2 + x y \left(-\frac{170}{9} x + 40 y + \frac{20}{3} z +1\right), \label{a4a4}
\end{equation}
has two real critical points $(0,0,0)$ and $\left( \frac{1}{64}, \frac{-1}{128}, \frac{-7}{192}\right)$ of class $A_4$,
 see \S~7.3 of the paper \cite{wall}. Let us perturb these critical points so that the resulting Morse function has only real critical points. The D-graph of this function can be split into two subgraphs of type $A_4$ in such a way  that all its edges connecting these subgraphs are oriented from one  subgraph to the other. Only one of the D-graphs \ref{54}--\ref{1233} admits such splitting.
The polynomial (\ref{a4a4}) is of class $\Xi_1$. Indeed, each its $A_4$ singularity can be perturbed so that the resulting function has no real critical points. According to \S \ref{morcom}, this is impossible for polynomials of type $\Xi_2$.

\subsubsection{Proof of Theorem \ref{refprop} for polynomials of type $\Xi_1$}
\label{prorefprop}

\begin{lemma} 
\label{lem44}
1. The set $\{E_6 + A_2\}$ of all third degree polynomials ${\mathbb R}^3 \to {\mathbb R}$ with non-discriminant principal homogeneous parts and two real critical points of classes $E_6$ and $A_2$ consists of exactly four smooth connected strata. 
Every polynomial with these properties can be reduced by an orientation preserving affine transformation of ${\mathbb R}^3$ and adding a constant
to one of the normal forms 
\begin{equation}
\frac{1}{3} x^3 \pm \frac{1}{3} z^3 - y^2 z + \varkappa z^2, \ \ \varkappa \neq 0,
\label{nfe6a2}
\end{equation}
where the sign $\pm $ at the monomial $\frac{1}{3} z^3$ and the sign of the coefficient $\varkappa$ characterize the stratum.
In particular, all these polynomials are of type $\Xi_1$. The $j$-invariant of their principal homogeneous parts is equal to 0.

2. The intersection of each of these four strata with the space 
$\Theta \cap {\mathbb R}^8$
of polynomials $($\ref{vers0}$)$ with real coefficients 
consists of three smooth surfaces diffeomorphic to ${\mathbb R}^2$. Namely, for any $\varkappa \neq 0$, all polynomials of the form $($\ref{vers0}$)$ that are affine equivalent to $\frac{1}{3} x^3 - \frac{1}{3} z^3 - y^2 z + \varkappa z^2$ are only the polynomial
\begin{equation}
\label{nexcl}
\frac{1}{3}(X^3+Y^3+z^3) + \beta (X+Y)^2,
\end{equation} 
where 
$\beta = \frac{\varkappa}{\sqrt[3]{16}}$, 
$X \equiv x- \beta,$ $Y \equiv y- \beta,$ 
and two other polynomials obtained from $($\ref{nexcl}$)$ by permutations of the coordinates.
Similarly, all polynomials of the form $($\ref{vers0}$)$ that are affine equivalent to
$\frac{1}{3} x^3 + \frac{1}{3} z^3 - y^2 z + \varkappa z^2$ are just the polynomial 
\begin{equation}
\frac{1}{3}(X^3+Y^3+Z^3) + 2X Y Z + \gamma (X-2Y+Z)^2,
\label{excl}
\end{equation}
where $\gamma = \frac{\varkappa}{\sqrt[3]{144}},$ 
$X \equiv x- \gamma,$ $Y \equiv y- 4\gamma,$ $Z \equiv z- \gamma,$ 
and two polynomials obtained from $($\ref{excl}$)$ by permutations of the coordinates.

Thus, for the parameters of the intersection surfaces we can in all cases take the coefficient $\varkappa \in (-\infty,0)$ or $(0, +\infty)$ and the constant that can be added to the polynomial $($\ref{nexcl}$)$ or $($\ref{excl}$).$
\end{lemma}

\noindent
{\it Proof.} 1. Let $f$ be such a polynomial. In coordinates centered at its point of class $E_6$, the quadratic part of its Taylor expansion has rank 1 and vanishes on a plane. The sign of this quadratic part is obviously an invariant of a component of polynomials of multi-class $\{E_6+A_2\}$, and the up-down involution of \S~\ref{invo} changes this sign. Therefore, it is sufficient to consider the case of non-negative quadratic parts. 

By \cite{LL}, Corollary 10.4.9, the line of zeros in ${\mathbb R}P^2$ of this quadratic part is tangent to the zero curve of the cubic part of $f$ at some inflection point of this curve. Choose the affine coordinates with the same origin so that this cubic part has the Weierstrass normal form (\ref{wei})
 with this inflection point equal to the line $\{y=z=0\}$. The quadratic part of $f$ is then equal to $\varkappa z^2$ with some $\varkappa>0$. Besides a critical point of multiplicity six at the origin, the resulting polynomial 
has two critical points in ${\mathbb C}^3$ which are the solutions of the system of equations $$y=0,  \quad 3x^2+a z^2=0, \quad 2a x + 3 b z + 2 \varkappa=0.$$ It is easy to calculate that these two critical points can only coincide if $a=0$. By a positive dilation of the coordinates, $f$ can then be reduced to the form (\ref{nfe6a2}).
 
Any polynomial (\ref{nfe6a2}) is mapped to itself by the reflection in the plane $\{y=0\}$. This reflection belongs to the non-unity component of the affine transformation group of ${\mathbb R}^3$. Therefore, any polynomial of the considered set $\{E_6 + A_2\}$ can be reduced to the form (\ref{nfe6a2}) by an orientation preserving affine transformation and adding a constant function. Thus, this set $\{E_6 + A_2\}$ consists of no more than four connected components characterized by combinations of signs of $\varkappa$ and $\pm$ in (\ref{nfe6a2}). The polynomials with different signs are topologically non-equivalent to each other (see Remark \ref{e67r}), so there are exactly four strata. 

By the well-known formula for the $j$-invariant (see e.g. \cite{BM}), the polynomials (\ref{wei}) have $a=0$ if and only if they define cubic curves with $j=0$. 

2. By another well-known formula (see \cite{BM}, formula (14)), the $j$-invariant of a real polynomial of the form (\ref{symnf}) is equal to 0 if and only if $A=0$ or $A=-2$. Thus, if $(x_0, y_0, z_0)$ are the coordinates of its critical point of class $E_6$, then in the local coordinates $X\equiv x-x_0, Y\equiv y-y_0, Z \equiv z-z_0$ this polynomial has one of the forms (\ref{nexcl}) or (\ref{excl}). 
The correspondence between two possible signs at the monomial $\frac{1}{3}z^3$ in (\ref{nfe6a2}) on the one hand and the values $0$ and $-2$ of the module $A$ in (\ref{symnf}) (i.e. the choice between the formulas (\ref{nexcl}) and (\ref{excl})) on the other hand follows from the comparison of the signatures of the quadratic parts of the critical points of class $A_2$ of the corresponding polynomials.

The expressions for the coefficients $\beta $ and $\gamma$ of the quadratic parts of the polynomials 
(\ref{nexcl}) and (\ref{excl}) by $\varkappa$ are obtained by comparing the critical values of these polynomials and (\ref{nfe6a2}) at the critical points of class $A_2$.

Finally, the coordinates $(x_0, y_0, z_0)$ of the critical points of class $E_6$ (involved in the expression of the coordinates $X, Y, Z$ through $x, y$ and $z$ in (\ref{nexcl}) and (\ref{excl})) follow from the requirement that the polynomials $($\ref{vers0}$)$ do not contain the monomials $x^2, y^2$ and $z^2$.
\hfill $\Box$ \smallskip

In the rest of the proof of Theorem \ref{refprop} we exploit the uniqueness of the stratum with both signs $+$ in (\ref{nfe6a2}).
\smallskip

Let $f!$ denote the polynomial (\ref{excl}) with $\beta=\sqrt[3]{16}$, so that its critical value at the critical point of class $A_2$ is equal to $\frac{4}{3}$.

Acting on the polynomials $f$ and $\tilde f$ from Theorem \ref{refprop} by appropriate orientation preserving affine transformations (which do not change their isotopy classes in the space of {\em generic} functions), we can assume that both $f$ and $\tilde f$ have the form (\ref{vers0}) with coefficients $(A, \lambda)$ from the parameter space $\Theta \cap {\mathbb R}^8$, where $\Theta \equiv ({\mathbb C}^1 \setminus \{1, e^{\pm 2\pi i/3}\}) \times {\mathbb C}^7$, see the proof of Theorem \ref{propmain}.

The systems of paths in ${\mathbb C}^1$ defining the vanishing cycles of $f$ and $\tilde f$ are homeomorphic to each other, so using the Lyashko--Looijenga covering over $B({\mathbb C}^1, 8)$ we can assume that $f$ and $\tilde f$ have the same collections of critical values and these systems of paths. 
Consider an arbitrary piecewise algebraic path $I:[0, 1] \to \Theta \cap {\mathbb R}^8$ consisting of real polynomials with principal parts of type $\Xi_1$, such that
\begin{enumerate}
\item $I(0) = f$,
\item $I(1) = f!$\ ,
\item the path $I$ has only finitely many transversal intersections with the variety of non-generic polynomials, 
\item for $t \in (0,1)$ sufficiently close to $1$ the points $I(t)$ are the Morse perturbations of the polynomial $f!$, which have only real critical points and realize the D-graph of Fig.~\ref{783}.
\end{enumerate}

Consider the path
$\Lambda \circ I: [0,1] \to \mbox{Sym}^8({\mathbb C}^1)$, i.e., the composition of $I$ and the Lyashko--Looijenga map. Let 
$\tilde I: [0,1] \to \Theta$, $\Lambda \circ \tilde I \equiv \Lambda \circ I$, be the lifting of this path $\Lambda \circ I$ to the space $\Theta$ with starting point $\tilde f$ (this lifting is uniquely defined according to the proof of Theorem \ref{propmain}). Since the virtual Morse functions of $f$ and $\tilde f$ are the same, the polynomials $I(t) $ and $\tilde I(t)$ undergo all the same elementary surgeries. In particular, this path $\tilde I$ lies in the real part $\Theta \cap {\mathbb R}^8$ of the space $\Theta$, and its endpoint is also a polynomial having a critical point of class $E_6$ with critical value $0$ and a point of class $A_2$, at which the polynomial is equivalent to $\xi^3 + \eta^2 - \zeta^2 + \frac{4}{3}$. By Lemma \ref{lem44}, this endpoint $\tilde I(1)$ can be reduced to $f!$ by an even permutation of the coordinates, which is an orientation preserving affine transformation of ${\mathbb R}^3$ and does not change the isotopy classes of the generic functions.
Applying this transformation to this endpoint and to the entire path $\tilde I$, we can assume that this endpoint coincides with the polynomial $f! \equiv I(1)$, and the polynomial $\tilde f$ is the starting point of this path. 

By the versality of the family (\ref{vers0}), a neighborhood of the point $f!$ in $\Theta$ is the space of a miniversal deformation of the multisingularity of class $(E_6 + A_2)$. Thus it has the structure of the direct product of the parameter spaces of the miniversal deformations of the singularities $E_6$ and $A_2$. The projections of this product onto the factors are defined by the {\em induction maps} (which are involved in the definition of the versality, see e.g. the map $\varphi$ in formula (1) of \S~8 in \cite{AVGZ1}) of this neighborhood (considered separately as the parameter space of a deformation of only the $E_6$ singularity and of only the $A_2$ singularity) onto these two parameter spaces of miniversal deformations.

Let $f_-$ and $\tilde f_-$ be near endpoints of the paths $I$ and $\tilde I$ that lie in this neighborhood of $f!$ and have equal sets of critical values. By the construction, the virtual Morse functions associated with them are the same. In particular, their two Coxeter--Dynkin graphs are the same and contain equal subgraphs with complementary sets of 6 and 2 vertices corresponding to the critical points of $f_-$ and $\tilde f_-$ obtained by decomposing two critical points of $f!$. Consider the projections of the points $f_-$ and $\tilde f_-$ to the parameter space of the standard miniversal deformation 
\begin{equation}
\label{vde6}
\tilde x^3 - \tilde y^4 + \tilde z^2 + \theta_1 + \theta_2 \tilde x + \theta_3 \tilde y + \theta_4 \tilde x \tilde y + \theta_5 \tilde y^2 + \theta_6 \tilde x \tilde y^2 
\end{equation}
of the $E_6$ singularity according to the structure of the direct product defined above. 
These projections belong to the same orbit of the group of automorphisms of the Lyashko--Looijenga covering over $B({\mathbb C}^1, 6)$ associated with this singularity. Indeed, the covering homotopy extends the map of one of these points to the other to an automorphism of the entire space of this covering, since the monodromy action of the fundamental group of the base on these points is controlled by their Coxeter-Dynkin diagrams, which are the same. This automorphism group is explicitly described in \cite{Liv}. It contains only two elements preserving the set of real functions: the identical map and the involution induced by the reflection $\tilde y \leftrightarrow -\tilde y$ of the argument space of the deformation (\ref{vde6}) and taking any point $(\theta_1, \theta_2, \theta_3, \theta_4, \theta_5, \theta_6)$ to $(\theta_1, \theta_2, -\theta_3, -\theta_4, \theta_5, \theta_6)$. 
The analogous real automorphism group for the $A_2$ singularity is trivial, so the projections of our two points $f_-$ and $\tilde f_-$ to the miniversal deformation of the $A_2$ singularity are the same. 

Thus, the point $\tilde f_-$ coincides either with $f_-$ or with the {\it unique} other point in the neighborhood of the point $f!$ which has the same projection as $f_-$ to the parameter space of the miniversal deformation of an $A_2$ singularity, and whose projection to the parameter space of the miniversal deformation of an $E_6$ singularity is obtained from that of $f_-$ by the unique non-trivial real automorphism of this deformation. But such a point is known to us: it is the image of $f_-$ under the reflection of the arguments in the plane $\{x=z\}$, so in this case $\tilde f_-(x, y, z) \equiv f_-(z, y, x)$.

In the first case (when $f_- =\tilde f_-$) we continue this equality along our paths and have $f = \tilde f$. In the second case, the identity 
$\tilde f_\tau (x, y, z) \equiv f_\tau (z, y, x)$ holds by the continuity over our paths for all their points, including the starting points $f$ and $\tilde f$, because the Lyashko--Looijenga covering commutes with the involution induced by the reflection $x \leftrightarrow z$. \hfill $\Box$

\subsubsection{Proof of Theorem \ref{achi}.} 
\label{prachi}
1. All virtual components mentioned in statement 1 of this theorem contain 
virtual Morse functions associated with polynomials which are invariant under reflections in some planes in ${\mathbb R}^3$, see formulas (\ref{twop})--(\ref{fourp}) and (\ref{reali77}). 

2. The D-graph of Fig.~\ref{1233} has no non-trivial symmetries
preserving the colors of the vertices, so all four critical points with Morse index 1 of any Morse polynomial with this D-graph can be ordered in an invariant way depending continuously of these polynomials. By Lemma \ref{flem} these four ordered points never lie in the same affine plane, and thus define an orientation of ${\mathbb R}^3$. The polynomials of the same isotopy class define equal orientations, but any polynomial and its mirror image define opposite orientations. \hfill $\Box$

\subsection{Polynomials with six real critical points}

\begin{proposition}
\label{prosx}
The set of all 1512 virtual Morse functions of type $\Xi_1$ with exactly six real critical points splits into five virtual components: one with cardinality 28 and passport $(1, 2, 2, 1)$, two with cardinality 112 and passports $(1, 3, 2, 0)$ and $(0, 2, 3, 1)$, and two with cardinality 630 and passport $(0, 3, 3, 0)$. Each of these virtual components is achiral and so is associated with a single isotopy class of Morse polynomials.
\end{proposition}

\noindent
{\it Proof.} A computer search analogous to the one described in the proof of Proposition \ref{pro14} proves that there are virtual components with six real critical points and values of Card invariant equal to 28, 112 and 630. 
Two virtual Morse functions 
\begin{equation}
\label{virtuM1}
\begin{array}{||cccccc|cc|}
\hline
 $-2$ & 0 & 0 & 1 & 1 & 0 & 1 & $-1$ \\
 0 & $-2$ & 0 & 0 & 1 & 0 & 0 & 0 \\
 0 & 0 & $-2$ & 0 & 1 & 1 & 1 & $-1$ \\
 1 & 0 & 0 & $-2$ & 0 & 0 & $-1$ & 0 \\
 1 & 1 & 1 & 0 & $-2$ & 0 & $-1$ & 1 \\
 0 & 0 & 1 & 0 & 0 & $-2$ & $-1$ & 0 \\
 1 & 0 & 1 & $-1$ & $-1$ & $-1$ & $ -2$ & 0 \\
 $-1$ & 0 & $-1$ & 0 & 1 & 0 & 0 & $-2$ \\
\hline
 0 & 0 & 0 & $-1$ & 1 & $-1$ & 0 & 0 \\
\hline
 o & o & o & e & e & e & & \\
\hline
\end{array} \qquad \qquad
\begin{array}{||cccccc|cc|}
\hline
 $-2$ & 0 & 0 & 0 & 0 & 1 & 0 & 0 \\
 0 & $-2$ & 0 & 1 & 0 & 0 & 0 & 0 \\
 0 & 0 & $-2$ & 1 & 1 & 1 & 1 & $-1$ \\
 0 & 1 & 1 & $-2$ & 0 & 0 & $-1$ & 0 \\
 0 & 0 & 1 & 0 & $-2$ & 0 & 0 & 1 \\
 1 & 0 & 1 & 0 & 0 & $-2$ & $-1$ & 0 \\
 0 & 0 & 1 & $-1$ & 0 & $-1$ & $ -2$ & $-1$ \\
 0 & 0 & $-1$ & 0 & 1 & 0 & $-1$ & $-2$ \\
\hline
 0 & 0 & 0 & 0 & $-1$ & 0 & $-1$ & $-1$ \\
\hline
 o & o & o & e & e & e & & \\
\hline
\end{array} 
\end{equation}
from this list have Card=630 (i.e., our restricted program quoted on p. \pageref{pro22} starting from any one of them reports that the corresponding virtual component consists of 630 elements).

By Proposition \ref{lemfun}, the first three critical values of any generic Morse polynomial associated with either of these two virtual Morse functions 
correspond to critical points whose Morse indices are not equal to 3 (because they are smaller than the indices of some other points). Being odd, these Morse indices are equal to 1. Similarly, the Morse indices of the other three real critical points are equal to 2. In particular, all such polynomials have passport $(0, 3, 3, 0)$.

 Let us prove that the virtual components represented by two virtual Morse functions (\ref{virtuM1}) are different.

At an additional request, our program finds all virtual Morse functions in the virtual component of the left virtual Morse function (\ref{virtuM1}) such that 
\begin{itemize}
\item[(i)] the number of non-zero intersection indices $\langle \Delta_i, \Delta_j \rangle$ with $i< j \leq 6$ (i.e., corresponding to real critical points) is at most 5;
\item[(ii)] the number of non-zero intersection indices $\langle \Delta_i, \Delta_j \rangle$ with $i \leq 6$ and $j \in \{7, 8\}$ (i.e., with $\Delta_i$ corresponding to a real critical point and $\Delta_j$ to non-real one) is at most 5;
\item[(iii)] the number of negative critical values is at most 4.
\end{itemize}

It turns out that there are exactly six virtual Morse functions in this component that satisfy all of these conditions, and for all of them the number of negative critical values is equal to 4. The right
virtual function (\ref{virtuM1}) satisfies the above conditions (i) and (ii) but has no negative critical values. Therefore, these two virtual Morse functions cannot belong to the same virtual component, and any two real polynomials associated with them are separated by the set-valued invariant.

A virtual component with Card=112 contains the virtual Morse function (\ref{virtuM2}) (left).
\begin{equation}
\label{virtuM2}
\begin{array}{|cccccc||cc|}
\hline
 $-2$ & 0 & 0 & 1 & 1 & $-1$ & $-1$ & 0 \\
 0 & $-2$ & 1 & 0 & 1 & $-1$ & $-1$ & 0 \\
 0 & 1 & $-2$ & 0 & 0 & 1 & 0 & 0 \\
 1 & 0 & 0 & $-2$ & 0 & 1 & 0 & 0 \\
 1 & 1 & 0 & 0 & $-2$ & 1 & 1 & $-1$ \\
 $-1$ & $-1$ & 1 & 1 & 1 & $-2$ & 0 & 0 \\
 $-1$ & $-1$ & 0 & 0 & 1 & 0 & $ -2$ & 0 \\
 0 & 0 & 0 & 0 & $-1$ & 0 & 0 & $-2$ \\
\hline
 1 & 1 & $-1$ & $-1$ & $-1$ & $2$ & 0 & 0 \\
\hline
 o & o & e & e & e & o & & \\
\hline
\end{array} 
\qquad \qquad
\begin{array}{|cccccc||cc|}
\hline
 $-2$ & 1 & 1 & $-1$ & $-1$ & 1 & 0 & 0 \\
 1 & $-2$ & 0 & 1 & 1 & $-1$ & $-1$ & 0 \\
 1 & 0 & $-2$ & 1 & 1 & $- 1$ & 0 & 1 \\
 $-1$ & 1 & 1 & $-2$ & 0 & 1 & 0 & 0 \\
 $-1$ & 1 & 1 & 0 & $-2$ & 1 & 1 & $-1$ \\
 1 & $-1$ & $-1$ & 1 & 1 & $-2$ & 0 & 0 \\
 0 & $-1$ & 0 & 0 & 1 & 0 & $ -2$ & 0 \\
 0 & 0 & 1 & 0 & $-1$ & 0 & 0 & $-2$ \\
\hline
 $-1$ & 1 & 1 & $-1$ & $-1$ & 2 & 0 & 0 \\
\hline
 e & o & o & e & e & o & & \\
\hline
\end{array}
\end{equation}
By Proposition \ref{lemfun}, its passport is equal to $(0, 2, 3, 1)$. Therefore, its up-down analog has passport $(1, 3, 2, 0)$, in particular
belongs to a different virtual component (but has the same Card invariant value). 

A virtual component with Card=28 contains the virtual Morse function (\ref{virtuM2}) (right). By Proposition \ref{lemfun}, its passport is equal to $(1, 2, 2, 1)$. 

The sum of the cardinalities of the five virtual components obtained is equal to the total number of virtual Morse functions of type $\Xi_1$ with six real critical points, therefore there are no other virtual components. By Corollary \ref{corud}, two virtual components 
with Card=630 are not invariant under the up-down involution of \S~\ref{invo}, therefore this involution maps them into each other.
\hfill $\Box$ \smallskip

Let us realize these five components by real polynomials.
The polynomial
\begin{equation}
\label{sxone}
x^3- 3 x + y^3 - 3 y + z^3- 3 \varepsilon (1+x+y)z
\end{equation}
with small $\varepsilon > 0$ has passport $(1, 3, 2, 0)$, its opposite polynomial has passport $(0, 2, 3, 1)$,
the polynomial
\begin{equation}
\label{sxthree}
x^3- 3 x + y^3 - 3 y + (1+x+y)z^2 + \varepsilon z^3
\end{equation}
with $|\varepsilon| \in (0, \frac{2}{3\sqrt{3}})$ has passport $(1, 2, 2, 1)$. By Proposition \ref{prosx}, the Card invariants of these three polynomials are equal respectively to 112, 112 and 28. All these polynomials are self-symmetric with respect to the plane $\{x=y\}$ and therefore their isotopy classes are achiral. 

Furthermore, we can perturb two critical points of the polynomial (\ref{e66}) in such a way that the point of class $A_2$ splits into two non-real critical points, and the point of class $E_6$ splits as in \S~ \ref{rea783} into three real Morse points of index 1 and three points of index 2 (thus providing the passport $(0, 3, 3, 0)$). This can be achieved by adding the term $\frac{\varepsilon^4}{40}x z^2$ to (\ref{reali77}). The resulting polynomial 
\begin{equation}
\label{reali77a}
z^2 + \frac{1}{3}x^3 + \frac{1}{3}z^3 - y^2 z - \varepsilon x y^2 - \frac{\varepsilon^4}{12}x - \frac{\varepsilon^3}{2}y^2 +\frac{\varepsilon^4}{40}x z^2
\end{equation}
is invariant under the reflection in the plane $\{y=0\}$, so its isotopy class is also achiral.
\hfill $\Box$

\subsection{Polynomials with four critical points}

\begin{proposition}
The set of all 515 virtual Morse functions of type $\Xi_1$ with exactly four real critical points splits into four virtual components: one with cardinality 30, two with cardinality 60, and one with cardinality 365. Each of these virtual components is achiral and associated with a single isotopy class of Morse polynomials.
\end{proposition}

\noindent
{\it Proof.} As before, our program finds three virtual components with Card equal to 30, 60 and 365. A virtual component with Card=60 contains the virtual Morse function (\ref{virtuM24}) (left).
\begin{equation}
\label{virtuM24}
\begin{array}{|cccc||cccc|}
\hline
 $-2$ & 1 & 1 & $-1$ & $-1$ & 0 & $-1$ & $-1$ \\
 1 & $-2$ & 0 & 1 & 0 & 0 & 0 & 0 \\
 1 & 0 & $-2$ & 1 & 0 & $-1$ & 0 & 1 \\
 $-1$ & 1 & 1 & $-2$ & 0 & 0 & 0 & 0 \\
 $-1$ & 0 & 0 & 0 & $-2$ & $-1$ & 0 & 0 \\
 0 & 0 & $-1$ & 0 & $-1$ & $-2$ & 0 & 0 \\
 $-1$ & 0 & 0 & 0 & 0 & 0 & $ -2$ & $-1$ \\
 $-1$ & 0 & 1 & 0 & 0 & 0 & $-1$ & $-2$ \\
\hline
 1 & $-1$ & $-1$ & 2 & 0 & 0 & 0 & 0 \\
\hline
 o & e & e & o & & & & \\
\hline
\end{array} 
\qquad \qquad
\begin{array}{|cccc||cccc|}
\hline
 $-2$ & 1 & 0 & 1 & 0 & $-1$ & $-1$ & 0 \\
 1 & $-2$ & 0 & 0 & 0 & 0 & 0 & 0 \\
 0 & 0 & $-2$ & 1 & 0 & 0 & $-1$ & 1 \\
 1 & 0 & 1 & $-2$ & $-1$ & 1 & 1 & $-1$ \\
 0 & 0 & 0 & $-1$ & $-2$ & 1 & 0 & 0 \\
 $-1$ & 0 & 0 & 1 & 1 & $-2$ & 0 & $-1$ \\
 $-1$ & 0 & $-1$ & 1 & 0 & 0 & $ -2$ & 1 \\
 0 & 0 & 1 & $-1$ & 0 & $-1$ & 1 & $-2$ \\
\hline
 0 & 0 & 1 & 0 & 0 & 1 & 0 & 1 \\
\hline
 o & e & o & e & & & & \\
\hline
\end{array} 
\end{equation}
By Proposition \ref{lemfun}, all generic real polynomials associated with this virtual Morse function have passport $(0, 1, 2, 1)$, so the up-down involution takes them to polynomials with passport $(1, 2, 1, 0)$ that have a different set-valued invariant but the same value of Card invariant. Since $30 + 60 \times 2 + 365 = 515,$ we have no additional virtual components.
\hfill $\Box$ \smallskip

Let us realize these components. The polynomial 
\begin{equation}
 x^3+y^3+z^3-3 z + 3 \varepsilon (x+y)z
\label{eee1}
\end{equation}
with sufficiently small $\varepsilon > 0$ has passport $(0, 1, 2, 1)$, minus this polynomial has passport $(1, 2, 1, 0)$, 
the polynomial 
\begin{equation}
\label{zh1}
x^3-3x +y^3-3y +z^3-3 \varepsilon z (1-(x+y)^2)
\end{equation}
has passport $(0, 2, 2, 0)$, and the polynomial
\begin{equation}
\label{zh2}
x^3-3x +y^3-3y +z^3+3 \varepsilon z (1-(x+y)^2)
\end{equation}
has passport $(1, 1, 1, 1)$.

All these polynomials are invariant under the reflection in the plane $\{x=y\}$, therefore they are achiral. 
The virtual component with Card=365 contains the virtual Morse function
(\ref{virtuM24}) (right), so by Proposition \ref{lemfun} its passport is equal to $(0, 2, 2, 0)$ and this component is realized by the polynomial (\ref{zh1}), while the component with Card=30 has passport $(1, 1, 1, 1)$ and is realized by the polynomial (\ref{zh2}).

\subsection{Polynomials with two critical points}

\begin{proposition}
The set of all 390 virtual Morse functions of type $\Xi_1$ with exactly two real critical points splits into three virtual components: two with cardinality 102 and one with cardinality 186. Each of these virtual components is associated with a single isotopy class of Morse polynomials.
\end{proposition}

\noindent
{\it Proof.} As before, this follows from the computer calculation. \hfill $\Box$ \medskip

The passport of the polynomial 
\begin{equation}
 x^3-3x +y^3-3y +z^3-3 \varepsilon z (x+y-1) 
\label{eee11}
\end{equation}
is equal to $(1, 1, 0, 0)$,
the passport of minus this polynomial is equal to $(0, 0, 1, 1)$,
and the passport of the polynomial  
\begin{equation}
\frac{1}{3}x^3 + \frac{1}{3}z^3 -y^2 z + z^2 +\varepsilon x (1+2z)
\label{eee14}
\end{equation}
with sufficiently small $\varepsilon >0$
is equal to $(0, 1, 1, 0)$. In the last case, adding the term $\varepsilon x(1+2z)$ to  polynomial (\ref{e66}) moves all of its critical points that are close to the origin into the non-real domain. It also splits the critical point of class $A_2$ into two real points with Morse indices 1 and 2. 

All these polynomials are invariant under certain reflections, so their isotopy classes are achiral.

\subsection{Polynomials without real critical points}

\begin{proposition}
There is only one virtual component of virtual Morse functions of type $\Xi_1$ without real critical points. This consists of 297 virtual Morse functions. 
\end{proposition}

\noindent
{\it Proof:} direct display of the program. \hfill $\Box$
 \smallskip

The corresponding real isotopy class is achiral and is represented by the polynomial (\ref{s}), $\varepsilon >0$.

\section{Polynomials of type $\Xi_2$}
\label{p82}

\subsection{Polynomials with eight critical points}

\begin{proposition}
\label{only2}
The set of all 6372 virtual Morse functions of class $\Xi_2$ with eight real critical points consists of exactly eight virtual components: 
\begin{itemize}
\item one of cardinality 324 with the D-graph shown in Fig.~\ref{324},

\item one of cardinality 504 with the D-graph shown in Fig.~\ref{504},

\item one of cardinality 945 with the D-graph shown in Fig.~\ref{945},

\item one of cardinality 1413 with the D-graph shown in Fig.~\ref{1413},
\end{itemize}
and four other virtual components obtained from these by the up-down involution, in particular with the same four cardinalities. 
\end{proposition}

\begin{figure}
\unitlength 0.7mm
\begin{center}
\begin{picture}(123,53)
\put(0,5){\circle{1.9}} 
\put(40,5){\circle{1.9}}
\put(80,5){\circle{1.9}}
\put(120,5){\circle{1.9}}
\put(0,45){\circle*{1.9}}
\put(40,45){\circle*{1.9}}
\put(80,45){\circle*{1.9}}
\put(120,45){\circle*{1.9}}
\put(0,6){\vector(0,1){38}}
\put(40,36){\vector(0,1){7}}
\put(40,34){\line(0,-1){8}}
\put(40,24){\line(0,-1){8}}
\put(40,14){\line(0,-1){8}}
\put(120,36){\vector(0,1){7}}
\put(120,34){\line(0,-1){8}}
\put(120,24){\line(0,-1){8}}
\put(120,14){\line(0,-1){8}}
\put(79.6,14){\vector(0,-1){7}}
\put(80.4,14){\vector(0,-1){7}}
\put(79.6,16){\line(0,1){8}}
\put(80.4,16){\line(0,1){8}}
\put(79.6,26){\line(0,1){8}}
\put(80.4,26){\line(0,1){8}}
\put(79.6,36){\line(0,1){8}}
\put(80.4,36){\line(0,1){8}}
\put(0.8,5.4){\vector(2,1){78}}
\put(40.7,5.7){\vector(1,1){38.5}}
\put(119.3,5.7){\vector(-1,1){38.5}}
\put(100.5,25.5){\vector(1,1){19}}
\put(99.5,24.5){\line(-1,-1){19}}
\put(52.9,32.1){\vector(-1,1){12.3}}
\bezier{100}(53.8,31.2)(55.8,29.2)(59.6,25.4)
\put(60.5,24.5){\line(1,-1){18.8}}
\put(81,45){\vector(1,0){38}}
\put(79,45){\vector(-1,0){38}}
\put(119,5){\vector(-1,0){38}}
\put(41,5){\vector(1,0){38}}
\bezier{600}(0.7,4.3)(60,-5)(78,4)
\put(78,4){\vector(2,1){1}}
\end{picture}
\end{center}
\caption{$D_5 + A_3$, Card = 32}
\label{324}
\end{figure}

\begin{figure}
\unitlength 0.7mm
\begin{center}
\begin{picture}(118,88)
\put(60,5){\circle*{1.9}} 
\put(0,45){\circle{1.9}} 
\put(40,45){\circle{1.9}} 
\put(80,45){\circle{1.9}} 
\put(0,85){\circle*{1.9}} 
\put(40,85){\circle*{1.9}} 
\put(80,85){\circle*{1.9}} 
\put(113.5,45){\circle{1.9}} 
\put(59.3,5.7){\vector(-3,2){58.1}}
\put(59.5,5.8){\vector(-1,2){19}}
\put(60.5,5.8){\vector(1,2){19}}
\put(0.7,45.7){\vector(1,1){38.5}}
\put(40.7,45.7){\vector(1,1){38.5}}
\put(60.8,5.6){\vector(4,3){51.7}}
\put(0,46){\vector(0,1){38.3}}
\put(19.5,65.5){\vector(-1,1){18.5}}
\put(20.5,64.5){\line(1,-1){18.5}}
\put(80,46){\vector(0,1){38.5}}
\put(59.5,65.5){\vector(-1,1){18.6}}
\put(60.5,64.5){\line(1,-1){18.5}}
\put(15,65){\vector(-3,4){14.5}}
\put(21,57){\line(3,-4){15}}
\put(42,29){\line(3,-4){16}}
\put(47,57){\vector(-1,4){6.7}}
\put(48.5,51){\line(1,-4){5}}
\put(55.5,23){\line(1,-4){3.5}}
\put(73,57){\vector(1,4){6.7}}
\put(71.5,51){\line(-1,-4){5}}
\put(65,25){\line(-1,-4){4.8}}
\end{picture}
\end{center}
\caption{Other side of $D_5 + A_3$, Card = 504}
\label{504}
\end{figure}

\begin{figure}
\unitlength 0.7mm
\begin{center}
\begin{picture}(123,48)
\put(0,5){\circle{1.9}} 
\put(40,5){\circle{1.9}}
\put(120,5){\circle{1.9}}
\put(120,45){\circle*{1.9}}
\put(0,45){\circle*{1.9}}
\put(40,45){\circle*{1.9}}
\put(80,45){\circle*{1.9}}
\put(80,5){\circle{1.9}} 
\put(0,6){\vector(0,1){38.5}}
\put(0.8,5.6){\vector(2,1){78.5}}
\put(40,26){\vector(0,1){18}}
\put(40,24.6){\line(0,-1){18.5}}
\put(40.7,5.7){\vector(1,1){38.5}}
\put(119.3,5.7){\vector(-1,1){38.5}}
\put(120,6){\vector(0,1){38.5}}
\put(41,5){\vector(1,0){38}}
\put(119,5){\vector(-1,0){38}}
\bezier{600}(0.7,4.3)(50,-5)(79,4)
\put(79,4){\vector(2,1){0.4}}
\put(80.4,14){\vector(0,-1){8}}
\put(79.6,14){\vector(0,-1){8}}
\put(80.4,16){\line(0,1){8}}
\put(79.6,16){\line(0,1){8}}
\put(80.4,26){\line(0,1){8}}
\put(79.6,26){\line(0,1){8}}
\put(80.4,36){\line(0,1){8}}
\put(79.6,36){\line(0,1){8}}
\end{picture}
\end{center}
\caption{Card = 945}
\label{945}
\end{figure}

\begin{figure}
\unitlength 0.7mm
\begin{center}
\begin{picture}(123,53)
\put(0,5){\circle{1.9}} 
\put(40,5){\circle{1.9}}
\put(80,5){\circle{1.9}}
\put(120,5){\circle{1.9}}
\put(0,45){\circle*{1.9}}
\put(40,45){\circle*{1.9}}
\put(80,45){\circle*{1.9}}
\put(120,45){\circle*{1.9}}
\put(0,6){\vector(0,1){38.5}}
\put(0.7,5.7){\vector(1,1){38.5}}
\put(1,5){\vector(1,0){38.50}}
\put(41,45){\vector(1,0){38.5}}
\put(120,6){\vector(0,1){38.5}}
\put(40.8,5.6){\vector(2,1){78.50}}
\bezier{600}(40.7,45.5)(100,55)(119,46)
\put(119,45.5){\vector(3,-2){0.5}}
\put(40.7,5.7){\vector(1,1){38.50}}
\put(80,26){\vector(0,1){18}}
\put(80,24.6){\line(0,-1){18.5}}
\put(40.4,36){\line(0,1){8}}
\put(39.6,36){\line(0,1){8}}
\put(40.4,34){\line(0,-1){8}}
\put(39.6,34){\line(0,-1){8}}
\put(40.4,24){\line(0,-1){8}}
\put(39.6,24){\line(0,-1){8}}
\put(40.4,14){\vector(0,-1){8}}
\put(39.6,14){\vector(0,-1){8}}
\end{picture}
\end{center}
\caption{$A_5 + A_3$, Card = 1413}
\label{1413}
\end{figure}

\noindent
{\it Proof.} As in the proof of Proposition \ref{pro14}, our program finds four virtual components with these cardinalities and D-graphs. None of these D-graphs is invariant under the up-down involution described in Proposition \ref{invofd}, so there are four other virtual components with the same values of Card. The sum of these values over all eight components is equal to the number 6372 of all virtual Morse functions of type $\Xi_2$ with eight real critical points, so there are no additional virtual components. \hfill $\Box$

\begin{proposition}
\label{protrue2}
All white vertices of four D-graphs of Figs.~\ref{324}--\ref{1413} correspond to critical points of Morse index 1 of arbitrary Morse polynomials with these D-graphs. All their black vertices correspond to critical points of Morse index 2, except only for the lowest vertex in Fig.~\ref{504}, which corresponds to a local minimum.
\end{proposition}

\noindent 
{\it Proof.} This follows almost immediately from Proposition \ref{lemfun}. \hfill $\Box$

\subsubsection{Realization of the D-graph of Fig.~\ref{324}} 
\label{d5a32}
For any $\lambda \neq 0$, consider the polynomial 
\begin{equation}
\label{d5a3}
x^3 - x z^2 - y^2 z + \lambda x^2 
\end{equation} 
of type $\Xi_2$. Its principal (lower) quasihomogeneous part at the origin, 
\begin{equation}
\label{o5}
\lambda x^2 -x z^2- y^2 z \equiv \lambda \left(x-\frac{z^2}{2 \lambda}\right)^2 -\frac{1}{4 \lambda}z^4 - y^2 z, 
\end{equation}
 is of class $D_5$, so by the techniques of \cite{AVGZ1}, \S~12, the whole of this function germ is of this class. For $\lambda>0$ this critical point can be split by a small perturbation into three real Morse critical points with Morse index $1$ and two critical points with index $2$: 
\begin{figure}
\unitlength 0.9mm
\begin{center}
\begin{picture}(40,54)
\bezier{200}(0,50)(20,33)(30,33)
\bezier{200}(0,30)(20,47)(30,47)
\bezier{150}(30,33)(40,33)(40,40)
\bezier{150}(40,40)(40,47)(30,47)
\put(3,28){\line(0,1){24}}
\put(13,40){\circle{1.5}}
\put(3,32.5){\circle{1.5}}
\put(3,47.5){\circle{1.5}}
\put(28,10){\circle*{1.5}}
\put(8,10){\circle*{1.5}}
\put(16,10){\circle{1.5}}
\put(2.3,2.4){\circle{1.5}}
\put(2.3,17.6){\circle{1.5}}
\put(16.5,10){\vector(1,0){11}}
\put(15.5,10){\vector(-1,0){7}}
\put(2.8,2.9){\vector(3,4){4.8}}
\put(2.8,17.1){\vector(3,-4){4.8}}
\put(5,39){\small $+$}
\put(27,39){\small $+$}
\end{picture} \qquad \qquad \qquad \qquad 
\begin{picture}(40,54)
\bezier{200}(10,50)(30,33)(40,33)
\bezier{200}(10,30)(30,47)(40,47)
\bezier{150}(40,33)(50,33)(50,40)
\bezier{150}(50,40)(50,47)(40,47)
\put(23,40){\circle{1.5}}
\put(40,33){\circle{1.5}}
\put(40,47){\circle{1.5}}
\put(40,28){\line(0,1){24}}
\put(45,10){\circle*{1.5}}
\put(40,3){\circle{1.5}}
\put(40,17){\circle{1.5}}
\put(31,10){\circle*{1.5}}
\put(22,10){\circle{1.5}}
\put(40.5,16.5){\vector(2,-3){4.2}}
\put(40.5,3.5){\vector(2,3){4.2}}
\put(31,10){\vector(4,-3){8.5}}
\put(31,10){\vector(4,3){8.5}}
\put(31,10){\vector(-1,0){8.5}}
\put(31,10){\line(1,0){4}}
\put(36,10){\line(1,0){4}}
\put(41,10){\vector(1,0){4}}
\put(32,39){\small $-$}
\put(43.5,39){\small $+$}
\end{picture}
\end{center}
\caption{Morsifications of $D_5$ singularities}
\label{DDs}
\end{figure}
see Fig.~\ref{DDs} (left) for the zero level set of the corresponding splitting of a $D_5$ singularity of a function in 
two variables. The intersection indices of the standard basic vanishing cycles associated with these critical points can be computed by the Gusein-Zade--A'Campo method \cite{AC}, \cite{GZ} and are represented by the standard Dynkin graph of $D_5$ singularities: see the lower left part of Fig.~\ref{DDs}. 

For $\lambda<0 $ this critical point of class $D_5$ splits analogously into three critical points of Morse index $2 $ and two critical points of index $1$. 

In addition, the polynomial (\ref{d5a3}) has a critical point at $\left(-\frac{2 \lambda}{3}, 0, 0\right)$. In the local coordinates $\tilde x \equiv x + \frac{2\lambda}{3},$ $ y $ and $z$ at this point it is equal to 
$$\frac{4\lambda^3}{27} - \lambda \tilde x^2 + \frac{2 \lambda}{3} z^2 + \tilde x^3 - \tilde x z^2 - y^2 z.$$ The main (lower) quasihomogeneous part of its non-constant part is $$-\lambda \tilde x^2 + \frac{2\lambda}{3} z^2 - y^2 z \equiv - \lambda \tilde x^2 + \frac{2\lambda}{3}\left(z-\frac{3}{4\lambda}y^2 \right)^2 - \frac{3}{8 \lambda} y^4,$$
so the germ of our polynomial at this point is locally equivalent to 
$$\frac{4 \lambda^3}{27} + \lambda \left(-\tilde x^2 + \tilde z^2 - \tilde y^4\right), $$ in particular it is of class $A_3$.
For $\lambda >0$ it can be split by a small perturbation of the polynomial (\ref{d5a3}) into two critical points with Morse index $2$ and one point with Morse index $1$; for $\lambda<0$ vice versa. 

Let us perturb both critical points of the polynomial (\ref{d5a3}) simultaneously in the manner described. For any choice of $\lambda \neq 0$ in (\ref{d5a3}), the passport of the resulting polynomial is equal to $(0, 4, 4, 0)$. By an additional small perturbation, we can make the resulting Morse polynomial generic. Any Morse polynomial constructed in this way for some value $\lambda \neq 0$ is transformed by the involution (\ref{invof}) into a polynomial with the opposite value of $\lambda$. 

\begin{proposition}
\label{iden21}
The D-graph of the Morse polynomial just constructed is shown in Fig.~\ref{324} if $\lambda > 0$, or is obtained from it by the involution described in Proposition \ref{invofd} if $\lambda < 0$.
\end{proposition}

\noindent
{\it Proof.} By Proposition \ref{only2}, the D-graph of this polynomial for $\lambda > 0$ must be one of the four D-graphs of Figs.~\ref{324}--\ref{1413}, or one of the four D-graphs obtained from them by the involution described in Proposition \ref{invofd}. By the construction of this polynomial, the set of vertices of this D-graph can be split into two subsets of cardinalities 5 and 3 in such a way that, by removing all edges with endpoints in different subsets, we obtain the canonical Coxeter-Dynkin graphs of classes $D_5$ and $A_3$, and all removed edges are directed from the graph of class $D_5$ to the graph of class $A_3$. The D-graph shown in Fig.~\ref{324} allows such a splitting, and the remaining seven D-graphs do not, so our polynomial realizes exactly this D-graph. For the same reasons, our polynomial with $\lambda<0$ realizes the D-graph obtained from Fig.~\ref{324} by the involution described in Proposition \ref{invofd}. 

\subsubsection{Realization of Fig.~\ref{504}} We can perturb the critical point of class $D_5$ of the polynomial (\ref{d5a3}) with $\lambda >0$ in a different way than in the previous subsubsection, see Fig.~\ref{DDs} (right). This perturbation splits it into three Morse critical points with Morse index 1, one point with index 0 and one point with index 2. The corresponding Coxeter-Dynkin graph is shown in the lower right part of Fig.~\ref{DDs}. 

If we also perturb the critical point of class $A_3$ as before, we obtain a polynomial with passport $(1, 4, 3, 0)$.

As in the proof of Proposition \ref{iden21}, the D-graph of this polynomial can only be the one shown in Fig.~\ref{504}, and the D-graph of the similar polynomial with $\lambda<0$ is obtained from Fig.~\ref{504} by the involution described in Proposition \ref{invofd}. 

\subsubsection{Realization of Fig.~\ref{945}}
Consider the polynomial
\begin{equation}
\label{e6a1a1}
x^3 - x z^2 - y^2 z + \gamma z^2, \ \ \gamma \neq 0, 
\end{equation}
of type $\Xi_2$.
It has a critical point of class $E_6$ at the origin, and additionally two Morse critical points with Morse indices 1 and 2. The critical values at the last two points are equal to each other and to $\gamma^3$. For $\gamma > 0$ this polynomial has the form $u^3 - v^4 + w^2$
in some local coordinates at the origin. This its critical point can be split into three Morse points with Morse index 1 and three points with Morse index 2 (see \S~\ref{rea783}). The intersection matrix of the vanishing cycles of the corresponding Morse polynomial is computed in \S~7 of \cite{Petr9} and indeed defines the D-graph of Fig.~\ref{945}. 

\subsubsection{Realization of Fig.~\ref{1413}}
\label{a5a3}
The polynomial 
\begin{equation}
\label{A5A3}
f= x^3 - 4x z^2 + \frac{80}{27} z^3 - y^2 z + z(x-z) 
\end{equation}
is of type $\Xi_2$, since the polynomial $t^3 - 4 t + \frac{80}{27}$
has three real roots.

The polynomial (\ref{A5A3}) has critical points $0$ and $ \left(\frac{1}{4}, 0, \frac{3}{8}\right)$. Its Taylor decomposition in the local coordinates $\tilde x \equiv x - \frac{1}{4}, y $ and $\tilde z \equiv z-\frac{3}{8}$ at the latter point is  
$$ - \frac{1}{64} + \frac{3}{4}\left(\tilde x - \frac{4}{3} \tilde z\right)^2 - \frac{3}{8} y^2 + \tilde x^3 - 4 \tilde x \tilde z^2 + \frac{80}{27} \tilde z^3 - y^2 \tilde z .$$ 

The substitution $X = \tilde x - \frac{4}{3}\tilde z$ reduces it to the form 
$$-\frac{1}{64} + \frac{3}{4} X^2 - \frac{3}{8} y^2 + X\left(X^2 + 4X \tilde z +\frac{4}{3} \tilde z^2 \right) - y^2 \tilde z. $$ The main quasihomogeneous part of this function germ (minus the constant $- \frac{1}{64}$) with weights $(2, 2, 1)$ is equal to $$ \frac{3}{4}X^2 -\frac{3}{8} y^2 + \frac{4}{3} X \tilde z^2 \equiv \frac{3}{4} \left(X + \frac{8}{9} \tilde z^2 \right)^2 - \frac{3}{8} y^2 - \frac{16}{27} \tilde z^4 ,$$
so it is of class $A_3$. A small perturbation can split this critical point into two Morse critical points with Morse index 2 and one with Morse index 1.

The multiplicity of the system of equations $d f =0$ at the origin (i.e., the Milnor number of $f$ at this point) is at least 5.
Indeed, this system contains the equation $\frac{\partial f}{\partial y} \equiv -2yz =0.$ Substituting its solution $z=0$ into the other two equations gives a system whose four solutions are all at the origin; the result of substituting the solution $y=0$ also has at least one solution there. 

On the other hand, the Milnor number of $f$ at the origin is at most 5, since the sum of the Milnor numbers of all critical points in ${\mathbb C}^3$ of any polynomial of type $\Xi_2$ is equal to eight. The signature of the second differential of this critical point is obviously equal to $(1,1)$. Therefore, this is a point of class $A_5$, and the critical points obtained by any its morsification can only have Morse indices 1 and 2. The number of all critical points with even Morse indices of any morsification of $f$ is equal to the number of critical points with odd indices, so any decomposition of this $A_5$ singularity into five Morse critical points has three points with Morse index 1 and two points with index 2.

\begin{proposition}
The D-graph of the Morse polynomial obtained is shown in Fig.~\ref{1413}.
\end{proposition}

\noindent
{\it Proof.} Of all the four D-graphs shown in Figs.~\ref{324}--\ref{1413} and their involutions described in Proposition~\ref{invofd}, only the D-graph of Fig.~\ref{1413}
 admits a splitting of the set of its vertices into the subsets of cardinalities $3$ and $5$ in such a way that the interior subgraphs of these subsets are of types $A_3$ and $A_5$, and all edges of the D-graph connecting these subsets are directed from the subgraph $A_3$ to $A_5$. \hfill $\Box$ 

\subsubsection{Proof of Theorem \ref{refprop} for the case $\Xi_2$.}
\label{prorefprop2}

\begin{lemma} 
\label{lem55}
1. The set of all third degree polynomials ${\mathbb R}^3 \to {\mathbb R}$ with non-dis\-cri\-mi\-nant principal homogeneous parts and three critical points of classes $E_6$, $A_1$, and $A_1$ with coinciding critical values at two $A_1$ points consists of exactly two smooth connected strata of codimension 6 in the space of all third degree polynomials. All these polynomials can be reduced to the normal form $($\ref{e6a1a1}$)$ with $\gamma \neq 0$ by orientation preserving affine transformations of ${\mathbb R}^3$ and by adding constants; the strata are characterized by the sign of $\gamma$. In particular, all such polynomials are of type $\Xi_2$. The $j$-invariant of the principal homogeneous parts of all these polynomials is equal to 1728.

2. The intersection of each of these two strata with the space $\Theta \cap {\mathbb R}^8$ of real polynomials $($\ref{vers0}$)$ with $A \neq 1$ is the union of three smooth surfaces diffeomorphic to ${\mathbb R}^2$. One of these surfaces consists of all polynomials
\begin{equation}
\label{cur2}
\frac{1}{3}(X^3 + Y^3 + Z^3) - (1+\sqrt{3})X Y Z + \varkappa (X + (1+\sqrt{3})Y + Z)^2 + c, \ 
\end{equation}
where the sign of the coefficient $\varkappa \neq 0$ matches the sign of the coefficient $\gamma$ in $($\ref{e6a1a1}$)$ that characterizes the stratum,
$X \equiv x - \varkappa$, $Y \equiv y- (4+2\sqrt{3})\varkappa $ and $Z \equiv z- \varkappa$. Two other surfaces are obtained from this one by permutations of the coordinates.

3. The critical value of the polynomial $($\ref{cur2}$)$ at the Morse critical points is a third degree monomial of the parameter $\varkappa$.
\end{lemma}

\noindent
{\it Proof.} 1. Let $f$ be a polynomial that satisfies these conditions on its critical points and values. Proceeding as in the proof of Lemma \ref{lem44}, we arrive at the search for polynomials of the form $x^3 + a x z^2 + b z^3 - y^2 z + \theta z^2 ,$ $\theta \neq 0$, each of which has, in addition to a critical point of class $E_6$ at the origin, two different real critical points with equal critical values. The last condition easily implies that $b=0$ and $a<0$. By stretching the coordinates, we can make $a=-1$ (perhaps changing the coefficient at $z^2$) and thus reduce the polynomial to the form (\ref{e6a1a1}). 

The $j$-invariant of the principal part of (\ref{e6a1a1}) is equal to 1728. Conversely, a polynomial in the Weierstrass normal form (\ref{wei}) has $j=1728$ if and only if $b=0, a \neq 0$. If $a>0$ then it is of type $\Xi_1$, and if $a<0$ then it can be reduced to the form \ $x^3 - x z^2 - y^2z$ \ by dilating the coordinates. 

2. The numbers $1+\sqrt{3}$ and $1-\sqrt{3}$ are the only two values of $A$ for which the $j$-invariant of the polynomial (\ref{symnf}) is equal to 1728. However, this polynomial with $A=1-\sqrt{3}$ is of type $\Xi_1$.
The rest of the proof repeats that of statement 2 of Lemma \ref{lem44}.

3. The multiplicative group ${\mathbb R}_+$ of positive numbers acts on the space $\Xi_2$: for any $t>0$, the corresponding transformation is given by the formula 
\begin{equation}
\label{dil}
 f_{A, \lambda}(x, y, z) \mapsto t^3 f_{A, \lambda} (t^{-1} x, t^{-1} y, t^{-1} z). \hspace{3cm} 
\end{equation}
This transformation preserves the principal homogeneous part of the polynomial, multiplies the coefficient $\varkappa$ by $t$ and multiplies all critical values by $t^3$. \hfill $\Box$

The rest of the proof of Theorem \ref{refprop} follows the scheme of \S~\ref{prorefprop}. Namely, denote by $f!$ the unique polynomial of the form (\ref{cur2}) with critical value 1 at two real Morse critical points. 

Let $f$ and $\tilde f$ be two generic polynomials of type $\Xi_2$ with equal virtual Morse functions. Arguing as in \S~\ref{prorefprop}, we see that
there exist generic polynomials $g$ and $\tilde g$ of the form (\ref{vers0}) such that a) $g$ is isotopic to $f$ and $\tilde g$ is isotopic to $\tilde f$ in the class of generic polynomials, b) $g$ and $\tilde g$ have equal sets of critical values and standard paths in ${\mathbb C}^1$ defining the vanishing cycles, and c) there are two paths connecting $g$ and $\tilde g$ to $f!$ in the space of type $\Xi_2$ Morse polynomials of the form (\ref{vers0}) such that the images of these two paths under the Lyashko--Looijenga map coincide and their near endpoints are polynomials with the D-graph of Fig.~\ref{945}. It follows from the structure of a neighborhood of the point $f!$ in $\Theta$ that these near endpoints either coincide with each other (and then, by continuity along our paths, so do the starting points $g$ and $\tilde g$), or they are reduced to each other by the unique involution of this neighborhood preserving the Lyashko--Looijenga map: this involution is induced by the reflection $(x, y, z) \leftrightarrow (z, y, x)$ of ${\mathbb R}^3$. In the latter case, our paths are mapped to each other by this involution, in particular their starting polynomials $g$ and $\tilde g$ are also mapped to each other by the reflection of arguments in the plane $\{x=z\}$. \hfill $\Box$ 

\begin{theorem}
\label{achi2}
All eight virtual components mentioned in Proposition \ref{only2} are achiral. In particular, there are exactly eight isotopy classes of degree three Morse polynomials ${\mathbb R}^3 \to {\mathbb R}$ with non-discriminant principal homogeneous parts of type $\Xi_2$ and eight real critical points.
\end{theorem}

\noindent
{\it Proof.} Consider the virtual component, whose D-graph is shown in Fig.~\ref{324}. By Lemma \ref{p111}, there is a virtual Morse function in it, whose critical values are ordered in ${\mathbb R}^1$ as follows:
{\small
$\begin{array}{cccc}
4 & 8 & 5 & 7 \\
1 & 2 & 6 & 3 
\end{array}$
} (in this picture we refer to the positions of the corresponding vertices of the D-graph of Fig.~\ref{324}). Let $f$ be an arbitrary generic polynomial of the form (\ref{vers0}) with only positive critical values, associated with this virtual Morse function. Vertices 2 and 3 are not connected by edges of our D-graph, nor are vertices 7 and 8. Therefore, there exists a path in the real subspace of the parameter space $\Theta$ of the deformation (\ref{vers0}), starting from $f$, along which only two surgeries of type s2 occur, permuting the critical values of the corresponding pairs of critical points, and all the remaining four critical values are fixed along this path. The virtual Morse function of the resulting generic polynomial $\tilde f$ coincides with that of the initial polynomial $f$, so by Theorem \ref{refprop} either $f$ and $\tilde f$ are isotopic to each other in the space of generic polynomials ${\mathbb R}^3 \to {\mathbb R}$, or their isotopy classes are mapped to each other by the reflection of arguments in an arbitrary plane in ${\mathbb R}^3$. 

However, these polynomials cannot be isotopic in this space. Indeed, 
each isotopy class of generic polynomials with four critical points of Morse index 1 defines an orientation of ${\mathbb R}^3$ in the following way. We number these points in ascending order of their critical values. By Lemma \ref{flem} these points are the vertices of a simplex in ${\mathbb R}^3$, and their order defines an orientation of this simplex and of the whole space. This orientation is constant along the isotopy classes of generic polynomials. Given a generic path in an arbitrary isotopy class of Morse polynomials, these canonical orientations of the generic polynomials at the beginning and at the end of this path are the same if and only if the number of intersection points of this path with the set of non-strictly Morse polynomials corresponding to surgeries of type s2, which transpose critical values of two critical points with Morse index 1, is even. But the path constructed above provides only one such transposition. 

Thus, the polynomials $f$ and $\tilde f$ belong to different isotopy classes of generic polynomials, and hence (by Theorem \ref{refprop}) these isotopy classes (as well as the isotopy classes of the Morse polynomials containing them) are mapped to each other by the reflections of arguments in planes. On the other hand, $f$ and $\tilde f$ are connected by our path in the space of Morse polynomials, so the latter isotopy classes of Morse polynomials actually coincide, and this single isotopy class is mapped to itself by the reflection. This proves the achirality of the virtual component characterized by the D-graph of Fig.~\ref{324}.

For the virtual component with the D-graph shown in Fig.~\ref{504}, we use the same argumentation with the order
{\small
$
\begin{array}{cccc}
6 & 7 & 8 & \\
2 & 3 & 4 & 5 \\
& 1 & & 
\end{array}
$
}
and the path transposing the critical values (2, 3) and (7, 8).

For Fig.~\ref{945} we use the order 
{\small
$\begin{array}{cccc}
4 & 5 & 7 & 6 \\
1 & 2 & 8 & 3 
\end{array}
$
} and transpositions $(2, 3)$ and $(5,6)$. 

For Fig.~\ref{1413} we use the order 
{\small
$\begin{array}{cccc}
2 & 3 & 7 & 8 \\
1 & 6 & 4 & 5 
\end{array}
$ 
}
and transpositions $(4,5)$ and $(7,8)$. \hfill $\Box$

\subsection{Polynomials with six critical points}
\label{sxxxx}

\begin{proposition}
The set of all 1897 virtual Morse functions of type $\Xi_2$ with exactly six real critical points splits into four virtual components: two of cardinality 140, one of cardinality 665 and one of cardinality 952. 
\end{proposition}

\noindent
{\it Proof.} Our program says that there are virtual components with six real critical points and cardinalities equal to 140, 665 and 952. 
The passport of all polynomials associated with a virtual Morse function with $\mbox{Card}=140$ can be computed using Proposition \ref{lemfun} and is equal to $(1, 3, 2, 0)$, therefore there is another virtual component with the same value of the Card invariant and passport $(0, 2, 3, 1)$. Since $ 140 \times 2 + 665 + 952 = 1897,$ there are no additional virtual components.
\hfill $\Box$ \smallskip

Let us realize these components. 
\label{rea62}
The critical point of class $D_5$ of the polynomial (\ref{d5a3}) with $\lambda >0$ can be decomposed by an arbitrarily small perturbation
into three Morse critical points with Morse index $1$ and either two critical points with index $2$ or two points with indices $0$ and $2$ (see two parts of Fig.~\ref{DDs}). If we also split the point of class $A_3$ of the same polynomial into one real Morse critical point (which necessarily has Morse index 2) and two non-real points, then we obtain (depending on the previous choice) polynomials with passports $(0, 3, 3, 0)$ and $(1, 3, 2, 0)$. For example, these perturbations of the polynomial (\ref{d5a3}) with $\lambda=\frac{1}{4}$ are realized by polynomials
\begin{equation}
\label{d5a3def}
x^3 - x z^2 - y^2 z + \frac{1}{4} x^2 - 2 \varepsilon z^3 + 3 \varepsilon^2z^2 + 8 \varepsilon^3 z -2 \varepsilon y^2
\end{equation} 
and
\begin{equation}
\label{d5a3def2}
x^3 - x z^2 - y^2 z + \frac{1}{4} x^2 - \frac{\varepsilon}{2} z^3 + 3 \varepsilon^2z^2 + \frac{7 \varepsilon^3}{2} z - \frac{\varepsilon}{2} y^2 
\end{equation} 
with sufficiently small $\varepsilon > 0$. 
The up-down involution (\ref{invof}) turns the second of these into a polynomial with passport $(0, 2, 3, 1)$.

The perturbations 
\begin{equation}
\label{A5A3aa}
f= x^3 - 4x z^2 + \frac{80}{27} z^3 - y^2 z + z(x-z)+4 \varepsilon^4 y^2 - 5 \varepsilon^2 x^2 
\end{equation}
of the polynomial (\ref{A5A3}) with sufficiently small $\varepsilon >0$
split its critical point of class $A_5$ 
as in \S~\ref{a5a3} (i.e., into three points of Morse index 1 and two points of index 2), and the point of class $A_3$ into a single real critical point of index 2 and two non-real points. The passport of the polynomial obtained is equal to $(0, 3, 3, 0)$. 

\begin{proposition}
Two polynomials with passport $(0, 3, 3, 0)$ constructed in the previous two paragraphs have different values of the invariant Card, in particular they belong to different isotopy classes of Morse polynomials.
\end{proposition}

\noindent
{\it Proof.} These polynomials are obtained by single Morse surgeries of type s1 from polynomials with eight real critical points, whose D-graphs are shown in Figs.~\ref{324} and \ref{1413}, respectively. The virtual Morse functions associated with these polynomials are derived using the standard formulas (implemented in our program) from the similar data of the initial polynomials. These virtual Morse functions are
$$ \left|
\begin{array}{cccccc||cc}
\hline
 $-2$ & 0 & 0 & 0 & 1 & 0 & $1$ & 0 \\
 0 & $-2$ & 0 & 1 & 1 & 1 & 1 & $-1$ \\
 0 & 0 & $-2$ & 0 & $1$ & 0 & $0$ & $-1$ \\
 0 & 1 & 0 & $-2$ & 0 & 0 & $0$ & $0$ \\
 1 & 1 & 1 & 0 & $-2$ & 0 & $-1$ & $1$ \\
 0 & 1 & 0 & 0 & 0 & $-2$ & $-1$ & $1$ \\
 $1$ & 1 & $0$ & $0$ & $-1$ & $-1$ & $-2$ & 0 \\
 0 & $-1$ & $-1$ & $0$ & 1 & $1$ & 0 & $-2$ \\
\hline
 1 & $-1$ & 1 & 0 & 0 & 0 & 0 & 0 \\
\hline
 o & o & o & e & e & e & & \\
\hline
\end{array} \right| \qquad \mbox{and} \qquad
\left|
\begin{array}{cccccc||cc}
\hline
 $-2$ & 0 & $-2$ & 0 & 1 & 1 & 1 & 1 \\
 0 & $-2$ & 0 & 0 & 1 & 0 & 0 & 0 \\
 $-2$ & 0 & $-2$ & 0 & 1 & 1 & 1 & 1 \\
 0 & 0 & 0 & $-2$ & 0 & $1$ & $0$ & $0$ \\
 1 & 1 & 1 & 0 & $-2$ & 0 & 0 & 0 \\
 1 & 0 & 1 & $1$ & 0 & $-2$ & $0$ & $0$ \\
 1 & 0 & 1 & $0$ & 0 & $0$ & $-2$ & $-1$ \\
 1 & 0 & 1 & $0$ & 0 & $0$ & $-1$ & $-2$ \\
\hline 
 $0$ & $1$ & 0 & $1$ & 0 & $0$ & $-1$ & $-1$ \\
\hline
 e & o & o & o & e & e & & \\
\hline
\end{array} \right| .
$$
If we substitute them as initial data into our program that enumerates the elements of virtual components, we get the answer that the first of them belongs to a virtual component of cardinality 665, and the second to a virtual component of cardinality 952. \hfill $\Box$
 \smallskip

All polynomials (\ref{d5a3def})--(\ref{A5A3aa}) and the up-down version of (\ref{d5a3def2}) are invariant under the reflection in the plane $\{y=0\}$, therefore the corresponding virtual components are achiral and are associated with single isotopy classes of Morse polynomials.

\subsection{Polynomials with four critical points}

\begin{proposition}
The set of all 650 virtual Morse functions of type $\Xi_2$ with exactly four real critical points splits into three virtual components: two of cardinality 90 and one of cardinality 470. 
\end{proposition}

\noindent
{\it Proof.} Our program says that there are virtual components with four real critical points and cardinalities equal to 90 and 470. The passport of the polynomials associated with a virtual Morse function with $\mbox{Card}=90$ can be computed using Proposition \ref{lemfun} and is $(1, 2, 1, 0)$. So there is another virtual component with the same value of the Card invariant and the passport $(0, 1, 2, 1)$. \hfill $\Box$

The perturbation 
\begin{equation}
\label{d5a3def4}
x^3 - x z^2 - y^2 z + \frac{1}{4} x^2 - \varepsilon^2 z + \varepsilon y^2 
\end{equation}
of the polynomial (\ref{d5a3}) 
with $\lambda = \frac{1}{4}$ and sufficiently small $\varepsilon >0$ splits the $A_3$ point of (\ref{d5a3}) into two critical points of Morse index 2 and one of index 1; it also splits the $D_5$ point into a single point of index 1 and four non-real points, so that the resulting polynomial has passport $(0, 2, 2, 0)$.

Next, consider the perturbation 
\begin{equation}
\label{e6a1a1mod}
x^3 - x z^2 + y^2 z + \frac{1}{4} z^2 + 2 \varepsilon^2 y^2 + \varepsilon x y^2 - \varepsilon^4 x
\end{equation}
of the polynomial (\ref{e6a1a1}) with $\lambda=\frac{1}{4}$
and sufficiently small $\varepsilon>0$. It does not change the type of the two Morse critical points of 
(\ref{e6a1a1}) distant from the origin (with Morse indices 1 and 2) and breaks the singularity of class $E_6$ at the origin into two critical points with Morse indices 0 and 1 and four non-real critical points. So we get a polynomial with passport $(1,2,1,0)$; multiplying it by $-1$ we obtain a polynomial with passport $(0, 1, 2, 1)$ and the same value of the Card invariant. 

All three polynomials constructed in this subsection are invariant under the reflection in the plane $\{y=0\}$ and thus represent achiral Morse isotopy classes and virtual components.

\subsection{Polynomials with two critical points}

\begin{proposition}
The set of virtual Morse functions of type $\Xi_2$ with exactly two real critical points consists of only one virtual component of cardinality 255.
\end{proposition}

\noindent
{\it Proof.} Direct display of the program. \hfill $\Box$
 \smallskip

To realize this component, take the starting polynomial (\ref{e6a1a1}) of the realization of Fig.~\ref{945} and perturb its point of class $E_6$ in such a way that all critical points in its neighborhood move into the complex domain. This can be achieved by the polynomial
\begin{equation}
\label{e6a1a1moc}
x^3 - x z^2 + y^2 z + \lambda z^2 + \varepsilon x, \ \ \varepsilon > 0 .
\end{equation}
This polynomial is also invariant under the reflection in the plane $\{y=0\}$, thus the virtual component is achiral and is associated with only one isotopy class of Morse polynomials.

\subsection{Polynomials without real critical points} There are no such polynomials of class $\Xi_2$, see \S~\ref{morcom}. This fact is also confirmed by our program, see Proposition \ref{procount1}.

\begin{corollary}
There are no Morse polynomials of type $\Xi_2$ which simultaneously have
both a local minimum and a local maximum.
\end{corollary}

Indeed, this is true for all polynomials of this type constructed above, and hence also for all polynomials from their isotopy classes of Morse polynomials. \hfill $\Box$

\section{Isotopy classification of strictly Morse polynomials with eight real critical points}
\label{strict}

\begin{theorem}
\label{111}
For any virtual component of type $\Xi_1$ or $\Xi_2$ with eight real critical points, the number of isotopy classes of {\em strictly} Morse polynomials with this value of the set-valued invariant is equal to $\frac{2}{9}$ of the value of its Card invariant. Namely, the set of these isotopy classes splits into pairs mapped to each other by the reflection with respect to any plane in ${\mathbb R}^3$; these pairs are in one-to-one correspondence with the isomorphism classes of extensions of the natural partial order of vertices of the D-graph of our virtual component to a total order, cf. Lemma~\ref{p111}.
\end{theorem}

The first statement  of Theorem \ref{mthm2} is a direct corollary of this theorem and Proposition \ref{procount1}. \medskip

 For any generic polynomial $f$, denote by $f^\uparrow$ an arbitrary polynomial of the form $f + c$, where $c$ is a constant such that all real critical values of $f+c$ are positive. 

\begin{lemma}
\label{lemel}
1. Generic polynomials $f$ and $\tilde f$ are isotopic to each other in the class of strictly Morse polynomials if and only if $f^\uparrow$ and $\tilde f^\uparrow$ are isotopic in the class of strictly Morse polynomials whose critical values at all real critical points are positive.

2. If two generic polynomials whose critical points are all real and whose critical values are all positive are isotopic in the class of strictly Morse polynomials, then the virtual Morse functions associated with them are the same.

3. A generic polynomial $f(x, y, z)$ of degree three with only real critical points cannot be isotopic to its mirror image $f(x, -y, z)$ in the class of strictly Morse polynomials.
\end{lemma}

\noindent
{\it Proof.} This lemma coincides with Lemma 3 of \cite{rigidplane} (up to replacing two variables with three)  and has the same proof.
\hfill $\Box$

 \medskip
\noindent
{\it Proof of Theorem \ref{111}.} 
Consider an arbitrary virtual component $\Omega$, all elements of which have eight real critical points. By statement 1 of Lemma \ref{lemel}, the number of isotopy classes of strictly Morse polynomials with set-valued invariant $\Omega$ is equal to the number of isotopy classes of such polynomials with only positive critical values. 

The number of virtual Morse functions in the virtual component $\Omega$ that have only positive critical values, is equal to $\frac{1}{9} \mbox{Card}(\Omega).$ By Corollary \ref{cormain}, they are all represented by real polynomials. By statement 2 of Lemma \ref{lemel}, each isotopy class of strictly Morse polynomials with all positive critical values is associated with only one such virtual Morse function. Conversely, by Theorem \ref{refprop} each such virtual Morse function is associated with one or two such components. By statement 3 of Lemma \ref{lemel} it cannot be associated with one component, so we multiply $\frac{1}{9} \mbox{Card}(\Omega)$ by 2. \hfill $\Box$
\medskip

\noindent
{\it Proof of the second statement of Theorem \ref{mthm2}.} Denote by  
$\G$ the connected group of orientation preserving affine transformations $G: {\mathbb R}^3 \to {\mathbb R}^3$. This group acts on the spaces $\Xi_1$ and $\Xi_2$ by mapping each polynomial $f$ to $f \circ G$. The action of $\G$ on the set of strictly Morse polynomials with eight critical points is free. Indeed, if $f$ is such a polynomial and $f = f \circ G$, $G \in \G$, then the transformation $G$ maps each critical point of $f$ to itself. These points affinely generate ${\mathbb R}^3$, so $G = \mbox{Id}$.  According to the versality property of the deformation (\ref{vers0}), the orbits of this action intersect the space of polynomials (\ref{vers0}) transversally. Therefore, each connected component $U$ of the set of strictly Morse polynomials with eight critical points is a fiber bundle with fiber $\G$, whose base is an arbitrary connected component of the intersection of the component $U$ and the space (\ref{vers0}). This intersection is diffeomorphic to ${\mathbb R}^8$. Indeed, the Lyashko--Looijenga map makes it a covering over the domain $\{x_1 <x_2 < \cdots < x_8\} \subset {\mathbb R}^8$. Therefore, component $U$ is a fiber bundle with a contractible base, and its fiber is homotopy equivalent to $SL(3, {\mathbb R})$. \hfill $\Box$

\section{Two remarks}

\subsection{Normalization of D-graphs}

Consider one of the D-graphs of Figs. \ref{54}--\ref{1233} and \ref{324}--\ref{1413}, whose vertices are marked by {\it integer} Morse indices of the corresponding critical points (determined in Propositions \ref{protrue} and \ref{protrue2}).

\begin{definition} \rm
An edge of a D-graph with vertices marked by integer Morse indices is called {\em normal} if

a) it is oriented from a vertex with a smaller Morse index to a vertex with a larger index,

b) it is dashed if the parities of Morse indices of the critical points corresponding to its ends are the same; it is solid if these parities are different. 

 Otherwise, this edge is called a {\em tunnel} edge. The {\em normalization} of a $D$-graph consists in removing all its tunnel edges.
\end{definition}

Analyzing the D-graphs of all polynomials of type $\Xi_1$ or $\Xi_2$ with eight real critical points, we see that

a) all D-graphs of polynomials with at least one point of local minimum or maximum are already normal,

b) some eight out of ten D-graphs of polynomials with only critical points with Morse indices 1 and 2 are split by the normalization into pairs of standard Coxeter--Dynkin graphs of some simple singularities; 

c) each of the two remaining D-graphs (one of which is shown in Fig.~\ref{945} and the other is its up-down version) splits into an isolated vertex and the extended Coxeter--Dynkin graph \ 
\unitlength 0.5 mm
\begin{picture}(38,14)
\put(0,5){\circle*{3}}
\put(10,5){\circle*{3}}
\put(20,5){\circle*{3}}
\put(26,-3){\circle*{3}}
\put(26,13){\circle*{3}}
\put(36,-3){\circle*{3}}
\put(36,13){\circle*{3}}
\put(0,5){\line(1,0){10}}
\put(10,5){\line(1,0){10}}
\put(20,5){\line(3,4){6}}
\put(20,5){\line(3,-4){6}}
\put(26,-3){\line(1,0){10}}
\put(26,13){\line(1,0){10}}
\end{picture} \ of type $\tilde E_6$.
 \smallskip

A very similar situation arises for D-graphs of degree four polynomials ${\mathbb R}^2 \to {\mathbb R}$ and for perturbations of $J_{10}$ singularities, see \cite{rigidplane} and \cite{Vj10}. In particular, the extended Coxeter-Dynkin graphs of types $\tilde E_7$ and $\tilde E_8$ appear there.

\subsection{D-graphs and Morse complexes}

For each of our D-graphs shown in Figs. \ref{54}--\ref{1233} and \ref{324}--\ref{1413} and marked by the integer Morse indices, consider the 
sequence of vector spaces and homomorphisms 
$$ C_3 \to C_2 \to C_1 \to C_0$$
over ${\mathbb Z}_2$, where the generators of $C_i$ correspond to the vertices marked by $i$, and the matrix coefficient connecting two such generators of neighboring groups $C_i$ and $C_{i-1}$ is equal to the parity of the number of segments drawn between these vertices in the D-graph. 

In all cases, the homomorphisms defined by these matrix coefficients satisfy the chain identity $\partial \circ \partial =0,$
all arising chain complexes for D-graphs of type $\Xi_1$ are acyclic, and all complexes for D-graphs of type $\Xi_2$ have homology groups equal to ${\mathbb Z}_2$ in dimensions 1 and 2, cf. \S~\ref{morcom}.

\section{Proof of Theorems \ref{adj1} and \ref{adj2}}
\label{pt4}

\noindent
{\it Proof of Theorem \ref{adj1}.} A. The polynomials of type $\Xi_1$ with pairs of real critical points of classes $D_4^- + D_4^-$, $E_6 + A_2$, and $A_4 + A_4$ are constructed in \S~\ref{p8eight}, see formulas (\ref{zero}), (\ref{e66}), and  (\ref{a4a4}).

Every polynomial $x^3 + y^3 + z^3 - 3\varepsilon z$, $z>0$, has two critical points of type $D_4^+$.
It is easy to check that every polynomial $$x^3 + x z^2 - y^2 z + \lambda x^2, \quad \lambda \neq 0, $$
has two critical points of classes $D_5$ and $A_3$.

B. The polynomials of classes $\Xi_1$ or $\Xi_2$ with multisingularities $D_4^- + D_4^+$, $D_4^- + A_4$ and $D_6^- + A_2$ are forbidden by the Euler characteristic arguments. Namely, by \S~\ref{morcom} the sum of the local indices of the gradient vector field at all critical points of any polynomial of class $\Xi_1$ or $\Xi_2$ is equal to 0. These local indices are equal to 2 or $-2$ for $D_{2k}^-$ singularities, depending on the quadratic parts of these singularities, and are equal to 0 for $D_{2k}^+$ and $A_{2k}$ singularities.

The $A_7$, $D_7$ and $E_7$ singularities do not appear in the decompositions of the $P_8$ singularities (even in the complex domain), which excludes the multisingularities $A_7 + A_1$, \ $D_7+ A_1$, and $E_7 + A_1$ for both cases $\Xi_1$ and $\Xi_2$.

Also, none of the D-graphs shown in Figs.~\ref{54}--\ref{1233} can be split into two subgraphs whose inner edges form canonical Coxeter--Dynkin graphs of one of the pairs of classes $A_6$ and $A_2$, or $A_5$ and $A_3$, or
$A_4$ and the canonical Dynkin graph \
\unitlength 0.3 mm
\begin{picture}(33,24)
\put(0,9){\circle*{3}}
\put(30,9){\circle*{3}}
\put(15,-6){\circle*{3}}
\put(15,24){\circle*{3}}
\put(0,9){\line(1,1){15}}
\put(0,9){\line(1,-1){15}}
\put(30,9){\line(-1,-1){15}}
\put(30,9){\line(-1,1){15}}
\put(0,9){\line(1,0){8}}
\put(11,9){\line(1,0){8}}
\put(30,9){\line(-1,0){8}}
\end{picture} \
of the $D_4^+$ singularity, or
$A_2$ and the canonical Dynkin graph \ 
\begin{picture}(63,24)
\put(0,9){\circle*{3}}
\put(15,9){\circle*{3}}
\put(30,9){\circle*{3}}
\put(60,9){\circle*{3}}
\put(45,-6){\circle*{3}}
\put(45,24){\circle*{3}}
\put(0,9){\line(1,0){30}}
\put(30,9){\line(1,1){15}}
\put(30,9){\line(1,-1){15}}
\put(60,9){\line(-1,-1){15}}
\put(60,9){\line(-1,1){15}}
\put(30,9){\line(1,0){8}}
\put(60,9){\line(-1,0){8}}
\put(41,9){\line(1,0){8}}
 \end{picture} \
of the $D_6^+$ singularity
in such a way that all edges of our D-graph connecting vertices of different subgraphs are directed from one of these subgraphs to the other. However, such splittings would necessarily exist if a polynomial of type $\Xi_1$ had a pair of critical points of these types. \hfill $\Box$
 \smallskip

\noindent
{\it Proof of Theorem \ref{adj2}.} A. Polynomials of type $\Xi_2$ with pairs of real critical points of classes $D_5 + A_3$ and $A_5 + A_3$ are constructed in \S\S~\ref{d5a32} and \ref{a5a3}. 

B. If a polynomial of type $\Xi_2$ has two critical points of class $D_4^+$, then each of these critical points can be removed by a small perturbation, and we get a polynomial without real critical points, which contradicts \S~\ref{morcom}. So such polynomials do not exist. The same argument forbids polynomials of type $\Xi_2$ with multisingularities $A_4 + A_4$, $A_4 + D_4^+$, $E_6 + A_2$, $A_6 + A_2$, and $D_6^+ + A_2$.

Furthermore, none of the D-graphs shown in Figs.~\ref{324}--\ref{1413} can be split into two ordered subgraphs, whose inner edges form canonical Coxeter-Dynkin graphs 
\unitlength 0.4 mm
\begin{picture}(28,12)
\put(3,4){\circle*{2}}
\put(15,4){\circle*{2}}
\put(25,-2){\circle*{2}}
\put(25,10){\circle*{2}}
\put(3,4){\line(1,0){12}}
\put(15,4){\line(5,3){10}}
\put(15,4){\line(5,-3){10}}
\end{picture}
of class $D_4$, in such a way that all edges of the initial D-graph connecting the vertices of different subgraphs are directed from the first subgraph to the second. This prohibits the polynomials of class $\Xi_2$ with two critical points of class $D_4^-$. 
The absence of classes $D_4^- + D_4^+$, $D_4^- + A_4$, $D_6^- + A_2$, $A_7 + A_1$, $D_7 + A_1$, and $E_7 + A_1$ for polynomials of any type $\Xi_1$ or $\Xi_2$ is already proved in the proof of Theorem \ref{adj1}. \hfill $\Box$



\section{A few problems}

{\bf 1.} List all the isotopy classes of Morse polynomials within the parameter space of the canonical miniversal deformation (\ref{vers0}) of $P_8$ singularities. These classes are the connected components of the intersections of this space with the isotopy classes computed above. In general, there can be more than one such component. For example, in the case of the isotopy class of type $\Xi_1$ with D-graph shown in Fig.~\ref{783}, the boundary of each such component  contains only one polynomial of the form (\ref{vers0}) with critical points of types $E_6$ and $A_2$ and corresponding critical values 0 and 1. According to Lemma \ref{lem44}, there are  three such polynomials. One has the form  (\ref{excl}), and the  other two are obtained from it by coordinate permutation. Therefore, there are three different such components whose boundaries contain these three polynomials.
\medskip

{\bf 2.} Study the topology of the isotopy classes of Morse polynomials. The homotopy types of the caustic complements 
of simple singularities of classes $A_k$, $D_k$, and $E_6$ were computed in \cite{sed}, \cite{sede}.  All of them are either contractible or homotopy equivalent to the circles. In \cite{jsing}, it was proven that the caustic complements of parabolic singularities of classes $P_8^2$ and $J_{10}^3$ have nontrivial second homology groups. 

The following are some sources of homology classes of isotopy classes. \smallskip

A. For each isotopy class with passport $(m_0, m_1, m_2, m_3)$, its map to the configuration space of collections of \ $m_0+m_1+m_2+m_3$ \ distinct points in ${\mathbb R}^3$ colored into four colors (or fewer if some of the $m_i$ are zero) is defined in the obvious way. Some cohomology classes can be induced from the cohomology group of this configuration space. \smallskip

B. In the case of isotopy classes with non-real critical points, a ${\mathbb Z}_2$-valued 1-cohomology class is defined by the number of $s4$ surgeries occurring along the loops inside these classes. Can this cocycle be extended to an integer one? Can this cocycle take non-zero values on the isotopy classes of strictly Morse polynomials? \smallskip

C. Two ${\mathbb Z}_2$-valued 1-cohomology classes are also defined by the numbers of $s2$ surgeries at which two critical values at the Morse points with the same Morse index (equal to 1 or 2) become equal. \smallskip

D. If an isotopy class of Morse polynomials is achiral and does not contain polynomials that are invariant under the reflections in planes of ${\mathbb R}^3$, then it is the base of a nontrivial two-fold covering.  Therefore, it admits a nontrivial 1-cohomology class with ${\mathbb Z}_2$ coefficients. \medskip

{\bf 3.} Study the subdivision of isotopy classes of Morse polynomials with non-real critical values by the {\em complex Maxwell set}, which consists of the polynomials having non-real critical points with real critical values. \medskip

{\bf 4.} The first version of our program listing the virtual functions of singularities was written in Fortran in 1984; see \cite{AGLV2}. Revising and rewriting it it in a more advanced language would be useful. 

The bottleneck of its computational complexity lies in the procedure of comparing new virtual functions with previously found ones. Clustering the formal graph by combining at one vertex all virtual functions that differ only in the number of negative critical values  is likely to reduce computation time for a singularity with Milnor number $\mu$ by about a factor of $\mu^2$.

}

\end{document}